\documentclass[12pt]{amsproc}
\usepackage{amscd}
\usepackage{amssymb}

\newtheorem{theorem}{Theorem}[section]
\newtheorem{corollary}[theorem]{Corollary}
\newtheorem{lemma}[theorem]{Lemma}
\newtheorem{proposition}[theorem]{Proposition}

\theoremstyle{definition}
\newtheorem{definition}[theorem]{Definition}

\newtheorem{asss}[theorem]{Assumptions}

\newtheorem{rem}[theorem]{Remark}

\theoremstyle{remark}
\newtheorem{note}{Note\!}

\numberwithin{equation}{section}

\newsavebox{\SmallMathBox}

\DeclareRobustCommand*{\nicefrac}[2]{\ifmmode\mathnicefrac{#1}
{ #2}%
  \else\textnicefrac{#1}{#2}\fi}
\newcommand*{\textnicefrac}[2]{\check@mathfonts%
\mbox{\raisebox{.5ex}{\fontsize\sf@size\z@\selectfont#1}\kern-.
1em%
/\kern-.1em\raisebox{- .25ex}{\fontsize\sf@size\z@\selectfont#2}
}}
\newcommand*{\mathnicefrac}[2]{%
  \mathchoice
    {\m@fr@c{\scriptstyle}{#1}{#2}}
    {\m@fr@c{\scriptstyle}{#1}{#2}}
    {\m@fr@c{\scriptscriptstyle}{#1}{#2}}
    {\m@fr@c{\scriptscriptstyle}{#1}{#2}}}

\def\bbb{{\boldsymbol{\beta}}}

\newcommand{\abs}[1]{\lvert#1\rvert}

\def\Ci{C^\infty}
\def\dpa{\partial}

\newcommand{\dd}[1]{\frac{d}{d#1}}

\def\fequal#1{\stackrel{#1}{=}}

\def\finto#1{\stackrel{#1}{\into}}

\def\ftoo#1{\stackrel{#1}{\too}}

\def\ifff{\Longleftrightarrow}
\def\iii{{\boldsymbol{i}}}
\def\into{\hookrightarrow}
\def\lla{\langle}

\def\noi{\noindent}
\newcommand{\norm}[1]{\lVert#1\rVert}

\def\ol{\overline}

\def\rra{\rangle}
\def\sqm1{\sqrt{-1}}
\def\sss{{\mathsf{S}}}

\def\tand{\mbox{\ \rm  and }}

\def\too{\longrightarrow}

\def\wt{\widetilde}

\def\={\cong}
\def\>{\supset}
\def\<{\subset}
\def\ii{^{-1}}
\def\12{\frac{1}{2}}
\def\0{^{\circ}}

\def\AA{{\mathbb A}}

\def\CC{{\mathbb C}}
\def\EE{{\mathbb E}}

\def\MM{{\mathbb M}}
\def\NN{{\mathbb N}}

\def\RR{{\mathbb R}}

\def\ZZ{{\mathbb Z}}

\def\Aa{{\mathcal A}}
\def\Bb{{\mathcal B}}
\def\Cc{{\mathcal C}}

\def\Ff{{\mathcal F}}

\def\Ii{{\mathcal I}}

\def\Ll{{\mathcal L}}

\def\Pp{{\mathcal P}}

\def\Uu{{\mathcal U}}

\def\GGG{{\mathfrak G}}

\def\C{\CC}
\renewcommand*{\d}{\delta}
\newcommand*{\D}{\Delta}
\def\e{\varepsilon}
\def\f{\varphi}
\def\F{\Phi}
\def\g{\gamma}
\def\G{\Gamma}

\def\la{\lambda}

\def\m{\mu}

\def\N{\NN}

\def\R{\RR}

\def\s{\sigma}
\def\Si{\Sigma}

\def\w{\omega}

\def\z{\zeta}
\def\Z{\ZZ}

\DeclareMathOperator{\codim}{codim}

\DeclareMathOperator{\dom}{dom}

\DeclareMathOperator{\Grass}{Grass}

\DeclareMathOperator{\Hom}{Hom}

\DeclareMathOperator{\id}{id}

\def\index{\mbox{\rm index\,}}

\DeclareMathOperator{\Mas}{Mas}
\DeclareMathOperator{\mmax}{max}

\DeclareMathOperator{\ran}{im}
\DeclareMathOperator{\range}{im}

\DeclareMathOperator{\sa}{sa}
\DeclareMathOperator{\SF}{sf}

\begin{document}
\setcounter{page}{1} \setcounter{tocdepth}{2}

\renewcommand*{\labelenumi}{%
   (\roman{enumi})}

\title[Weak Symplectic Functional Analysis]
{Weak Symplectic Functional Analysis\\and General Spectral Flow Formula}


\author{Bernhelm Booss--Bavnbek}
\address{Institut for matematik og fysik\\ Roskilde
University, 4000 Ros\-kilde, Denmark} \email{booss@mmf.ruc.dk}

\author{Chaofeng Zhu}
\address{Nankai Institute of Mathematics\\
Nankai University\\
Tianjin 300071, P. R. China}
\email{zhucf@nankai.edu.cn}

\begin{abstract}
We consider a continuous curve of self-adjoint Fredholm extensions
of a curve of closed symmetric
operators with fixed minimal domain $D_m$ and fixed  {\it intermediate}
domain $D_W$\/.  Our main example is a family of symmetric generalized
operators of Dirac
type on a compact manifold with boundary with varying well-posed
boundary conditions. Here $D_W$ is the first Sobolev space and
$D_m$ the subspace of sections with support in the interior. We
express the spectral flow of the operator curve by the Maslov index
of a corresponding curve of Fredholm pairs of Lagrangian
subspaces of the quotient Hilbert space $D_W/D_m$ which is
equipped with continuously varying {\it weak symplectic structures}
induced by the Green form.

In this paper, we specify the continuity conditions; define the
Maslov index in weak symplectic analysis;
discuss the required weak inner Unique Continuation Property; derive a General
Spectral Flow Formula; and check that the assumptions are natural
and all are satisfied in geometric and pseudo-differential context.

Applications are given to $L^2$ spectral flow formulae; to the splitting
of the spectral flow on partitioned manifolds; and to linear Hamiltonian systems.
\end{abstract}

\maketitle


\tableofcontents


\section*{Introduction}

In various branches of mathematics one is interested in
the calculation of the spectral flow of a continuous family
of closed densely defined (not necessarily bounded) self-adjoint Fredholm
operators in a fixed Hilbert space. Roughly speaking, the spectral flow counts the
net number of eigenvalues changing from the negative real half axis to
the non-negative one.
The definition goes back to a famous paper by M.~Atiyah, V.~Patodi, and
I.~Singer \cite{AtPaSi75},
and was made rigorous by J.~Phillips \cite{Ph96} for continuous paths of
bounded self-adjoint Fredholm operators,
by K.P. Wojciechowski \cite{Wo85} and C. Zhu and Y. Long \cite{ZhLo99} in various
non-self-adjoint cases, and by B. Booss-Bavnbek, M. Lesch,
and J. Phillips \cite{BoLePh01} in the unbounded self-adjoint case.

\subsection{History}\label{ss:history}

The first spectral flow formula was the classical Morse index
theorem (cf. M. Morse \cite{Mo}) for geodesics on Riemannian
manifolds. It was generalized by W. Ambrose \cite{Am} in 1961 to
more general boundary conditions, which allowed two end points of
the geodesics varying in two submanifolds of the manifolds. In
1976, J.J. Duistermaat \cite{Du} completely solved the problem of
calculating the Morse index for the one-dimensional variational
problems, where the positivity of the second order terms was
required. In 2000-2002, P. Piccione and D.V. Tausk
\cite{PiT1,PiT2} were able to prove the Morse index theorem for
semi-Riemannian manifolds for the same boundary conditions as in
\cite{Am}, and some non-degenerate conditions were needed. In
2001, the second author \cite{Zh01} was able to solve the general
problem for the calculation of the Morse index of index forms for
regular Lagrangian systems. As mentioned in \cite{Zh01}, it is
nontrivial, in the general case, to see that this problem is
equivalent to the spectral flow formula of the corresponding
second order operators. Our Theorem \ref{t:hamiltonian} solves
this puzzle.

There is also another line of studying the spectral flow formula for first order operators.
Let $\{A_s:\Ci(M;E)\to\Ci(M;E)\}_{s\in [0,1]}$ be a family of continuously varying formally
self-adjoint linear elliptic differential operators of first order
over a smooth compact Riemannian manifold $M$ with
boundary $\Si$, acting on sections of a Hermitian vector bundle $E$ over $M$. Fixing a unitary
bundle isomorphism between the original bundle and a product bundle in a collar neighbourhood $N$
of the boundary, the operators $A_s$ can be written in
the form
\begin{equation}\label{e:intro-dirac}
A_s|_N = J_{s,t}(\frac{\dpa}{\dpa t}+B_{s,t})
\end{equation}
with invertible skew-self-adjoint bundle isomorphisms $J_{s,t}$ and
first order elliptic differential operators $B_{s,t}$ on $\Si$.
Here $t$ denotes the inward normal coordinate in $N$.

In 1988, A. Floer \cite{Fl88} studied the case that $M$ is a finite
interval with periodic boundary condition.
Later in 1991, T. Yoshida \cite{Yo91} studied the case $\dim M=3$ where $\{A_s\}$ is a curve of
Dirac operators with invertible ends, $J_{s,t}=J_s$
are unitary operators, $B_{s,t}=B_s$ symmetric.
In 1995, L. Nicolaescu \cite{Ni95} generalized Yoshida's results to arbitrary $\dim M$.
In 2000, M. Daniel \cite{Da00} removed the nondegenerate conditions in
\cite{Ni95}.
In 1998-2001, the first author, jointly with K. Furutani and N. Otsuki \cite{BoFu99,BoFuOt01}
proved the case that $A_s$ differ by
zeroth order operators, and the boundary condition is fixed.
In December 2000, P. Kirk and M. Lesch \cite{KiLe00} proved the case $A_s$ of Dirac type,
$J_{s,t}$ is fixed unitary, and $B_{s,t}=B_s$ symmetric.
In this paper we shall only assume that $B_{s,0}^*-B_{s,0}$ is a bundle homomorphism and
that $\ker A_s|_{H_0^1(M;E)}=\{0\}$ (weak inner unique continuation property - UCP).

Such operators we shall call {\em symmetric generalized operators of Dirac type}.
For details see Definition \ref{d:dirac}.

\smallskip

It may be worth mentioning that there is a multitude of {\em other} formulae involving
spectral flow, e.g., as error term under cutting and pasting of the index (see the first author with
K.P. Wojciechowski \cite[Chapter 25]{BoWo93}) or under pasting of the eta-invariant
as in \cite{KiLe00}. Whereas these formulae typically relate the spectral flow of a family
on a closed manifold of dimension $n-1$ to index or eta-invariant of a single operator on a
manifold of dimension $n$, this paper solely addresses the opposite direction,
i.e., how to express the spectral flow of a family over a manifold of dimension $n$
by objects (here by the Maslov index) defined on a hypersurface of dimension $n-1$.

\subsection{Our Setting, Goals, and Difficulties to  Overcome}
For some applications, the results obtained in \cite{BoFu99,BoFuOt01} respectively in \cite{KiLe00}
are not sufficient:
\begin{enumerate}
\item E.g., let $J_{s,t}$ and $B_{s,t}$ vary with $t$, admitted
in \cite{BoFu99,BoFuOt01}, but excluded in \cite{KiLe00}.

\item Let $J_{s,t}$ be only invertible, but not necessarily unitary, admitted
in \cite{BoFu99,BoFuOt01}, but excluded in \cite{KiLe00}.

\item Let $J_{s,t}$ vary with $s$, excluded both in \cite{BoFu99,BoFuOt01} and in \cite{KiLe00}.

\item Let the boundary conditions vary, excluded in \cite{BoFu99,BoFuOt01} and
admitted in \cite{KiLe00}.

\end{enumerate}

Our goal in this paper is to generalize both the results of
\cite{BoFu99,BoFuOt01} and of \cite{KiLe00}, i.e., we shall give a method to calculate
the spectral flow for a continuous curve of formally self-adjoint elliptic differential
operators of first order defined on a smooth manifold with
boundary with continuously varying self-adjoint elliptic (i.e., well-posed) boundary
conditions given by pseudo-differential projections.

To do that, the routes of \cite{BoFu99,BoFuOt01} and \cite{KiLe00} are barred to us
because they rely on the concept of one fixed symplectic Hilbert space.

The solution we give here is working in a smaller symplectic space of {\it reduced boundary values},
namely the quotient of a
fixed intermediate domain $D_W$ (the first Sobolev space in our applications)
and the minimal domain $D_m$ (the sections with support in the interior in our applications).
So, in terms of our applications, we have to work with the function space $D_W/D_m=H^{1/2}(\Si)$,
in addition to working with the natural distribution subspace
$\bbb=D_{\mmax}/D_m$ of $H^{-1/2}(\Si)$ and the familiar $L^2(\Si)$.
This leads to various difficulties:

\begin{enumerate}
\item In difference to $\bbb$ and to $L^2(\Si)$, Green's form
induces only a {\it weak} symplectic structure on the Hilbert space $H^{1/2}(\Si)$ (see Remark
\ref{r:weak-strong}b below).

\item There is no longer a canonical {\it splitting} (see Lemma \ref{l:strong_symplectic} below).

\item A priori, Lagrangian subspaces of our weak symplectic space do no longer characterize
{\it self-adjoint extensions} of a given closed symmetric operator.

\item The basic operator of our analysis is no longer the maximal closed extension but the
extension to the {\it intermediate} domain which clearly yields a bounded operator in suitable setting,
but {\it not a closed} operator
in $L^2$ in our applications (see our Assumption \ref{a:reduced-space} below).

\item A priori, continuity of domains in $\bbb$ or in $L^2(\Si)$ does not guarantee {\it continuity} in
$H^{1/2}(\Si)$ (see, however, Theorem \ref{t:cd_continuity} below).

\item A priori, possible differences between the {\it value of the Maslov index}
of the same curves in $\bbb$,
$L^2(\Si)$ and $H^{1/2}(\Si)$ can not be excluded (see, however, the
regularity inclusion \eqref{e:regularity} below).

\end{enumerate}
Consequently, in Sections \ref{s:weak}, \ref{s:sa-extensions}
our analysis requires to specify quite a list of supplementary assumptions. Fortunately,
it turns out in Section \ref{s:geometric} that all these assumptions are naturally satisfied
when we restrict ourselves to
symmetric generalized operators of Dirac type on smooth compact manifolds with boundary
and boundary conditions given by pseudo-differential projections.

\medskip

\subsection{Main Results}\label{ss:main-results}
The following four theorems will be proved.

Let $\{A_s: H^1(M;E)\to L^2(M;E)\}_{s\in [0,1]}$ be a family of
symmetric generalized operators of Dirac type over a compact
smooth Riemannian manifold $M$ with boundary $\Si$, acting on
sections of a Hermitian vector bundle $E$ over $M$. We assume that
the family is continuous as a family of bounded operators. Let
$\{Q_s: L^2(\Si;E|_{\Si})\to L^2(\Si;E|_{\Si})\}_{s\in [0,1]}$
denote the induced family of (pseudo-differential) Calder{\'o}n
projections. The family is continuous by an argument of
\cite[Section 3]{BoLePh01} in combination with M. Lesch
\cite{Le04a}. Generalizing the celebrated {\em Cobordism Theorem}
to our case (in preparation, see the present authors \cite{BoZh04}
in combination with \cite{Le04a}), we obtain that $\ran Q_s =
\g(\ker A_s^*)\cap L^2(\Si;E|_{\Si})$\/ is a Lagrangian subspace
of $L^2(\Si;E|_{\Si})$ endowed with a symplectic form defined by
$-J_{s,0}$\,. Here $\g:\dom(A^*_s)\to H^{-1/2}(\Si;E|_{\Si})$ is
the trace map. Let
\[
\{P_s: L^2(\Si;E|_{\Si})\to L^2(\Si;E|_{\Si})\}_{s\in [0,1]}
\]
be a continuous family of
pseudo-differential projections such that
\[
(\ker P_s,\ran Q_s)\in \Ff\Ll^2\bigl(L^2(\Si;E|_{\Si})\bigr) \text{ for all $s\in [0,1]$},
\]
where $\Ff\Ll^2\bigl(L^2(\Si;E|_{\Si})\bigr)$ denotes the space of Fredholm pairs
of Lagrangian subspaces of the (strong) symplectic Hilbert space $L^2(\Si;E|_{\Si})$.
(For better reading, we fix the underlying manifold $M$ and the bundle $E$.
For a more general setting in the terminology of fibre bundles, involving
bundle isomorphisms $T_s$, see Section \ref{ss:gen-di-op}.)

\begin{theorem}\label{t:gsff1} In the set-up described above,
spectral flow of the family $\{A_{s,P_s}\}$ and Maslov index
of the family of pairs
$\{\ker P_s\/,\range Q_s\}$ are well-defined, and we have the general spectral
flow formula
\[
\SF\{A_{s,P_s}\} = -\Mas\{\ker P_s,\range Q_s\}.
\]
Here $A_{s,P_s}$ acts like $A_s$ in $L^2(M;E)$ with
\[
\dom(A_{s,P_s}) = D_s := \{x\in H^1(M;E)\mid P_s(x|_{\Si})=0\}.
\]

\end{theorem}

\smallskip

Our versions of the spectral flow and the Maslov index are defined
in the Appendix in Definition \ref{d:spectral_continuous}c, respectively
below in Section \ref{ss:maslov} in  Definition \ref{d:continuous_banach}b.

\medskip

We obtain a second theorem in a slightly modified setting.
Now we assume that the manifold $M=M^+\cup_{\Si} M^-$ is
a partitioned closed manifold with a hypersurface $\Si$.
We denote the restrictions of $A_s$ to the parts by $A_s^{\pm}$\,.
Note that we now have a pair of Calder{\'o}n projections
$(Q_s^+,Q_s^-)$ for each $s\in [0,1]$ with $\range Q_s^\pm$
Lagrangian subspaces in $L^2(\Si;E|_{\Si})$ with symplectic form
again defined by $-J_{s,0}$.

\begin{theorem}\label{t:split}
For the partitioned case we assume that the parts operator $A_s^{\pm}$ are
symmetric generalized operators of Dirac type. Then we have
\[
\SF\{A_{s}\} = \SF\{A^-_{s,I-Q^+_s}\} = -\Mas\{\range Q^-_s,\range Q^+_s\}.
\]
\end{theorem}

\medskip

Next we consider a first order linear Hamiltonian system on the interval
$[0,T]$.
Let $j_{s,t}$, $b_{s,t}\in C^0([0,1]\times [0,T],{\rm gl}(m,\CC))$ be two families
of matrices such that $j_{s,t}^*=-j_{s,t}$ are invertible,
$b_{s,t}^*=b_{s,t}$, and $j_{s,t}$ is $C^1$ in $t$ for fixed $s$. Then we have
a family of unbounded Fredholm operators $\{A_s\}_{s\in [0,1]}$ in $L^2([0,T],\CC^m)$ defined by
$$(A_s x)(t):=-j_{s,t}\frac{d}{dt}x(t)-b_{s,t}x(t)-\frac{1}{2}
\left(\frac{d}{dt}j_{s,t}\right)x(t)$$
for all $s\in [0,1]$ with $x\in H^1([0,T],\CC^m)$.
We define the symplectic structure $\omega_s$
on $\CC^{2m}=\CC^m\oplus\CC^m$ by
$$\omega_s((x_1,x_2),(y_1,y_2)) := -\lla j_{s,0}x_1,y_1\rra +\lla j_{s,T}x_2,y_2\rra,$$
where $x_1,y_1,x_2,y_2\in\CC^m$.
Let $W_s$ be a continuous family of
Lagrangian subspaces of $(\CC^{2m}, \omega_s)$. Then the family
$\{A_{s,W_s}\}$ forms a continuous curve of self-adjoint Fredholm operators.
Here $A_{s,W_s}$ denotes the restriction of $A_s$ on the space
$$H_{W_s} := \{x\in H^1([0,T],\CC^m)|(x(0),x(T))\in W_s\}.$$
Let $\Gamma_s(t)$
be a fundamental solution of the linear system $A_sx=0$. Then we have

\begin{theorem}\label{t:First_Order} In the set-up described above,
the graphs $\{{\GGG}(\Gamma_s(T))\}$ of $\{\Gamma_s(T )\}$ make a continuous family
of Lagrangian subspaces of the symplectic vector space
$(\CC^{2m}, \omega_s)$. Moreover we have
$$\SF\{A_{s,W_s}\}=\Mas\{{\GGG}(\Gamma_s(T)),W_s\}.$$
\end{theorem}

\medskip

Now we consider regular linear Lagrangian systems. Let $p,q,r\in
C([0,1]\times[0,T],{\rm gl}(m,\CC))$ be families of matrices such that $p$
is of class $C^1$, $p_s(t)=p_s(t)^*$, $r_s(t)=r_s(t)^*$, and $p_s(t)$ is
invertible for all $s\in[0,1]$ and $t\in[0,T]$. Then we have a family of
formally self-adjoint linear differential operators of second order
$\{L_s\}_{s\in [0,1]}$ defined by
\begin{equation}\label{e:l-s}
L_s := -\frac{d}{dt}\left(p_s\frac{d}{dt}+q_s\right)+q_s^*\frac{d}{dt}+r_s.
\end{equation}
Let $J := \left(\begin{array}{ll}0&-I_m\\ I_m&0\end{array}\right)$, where
$I_m$ denotes the identity matrix on $\R^m$. When there is no confusion we
will omit the subindex of the identity matrices. Let

\begin{equation} \label{eq01m1.1}
b_s(t)=\left(\begin{array}{ll}
p_s^{-1}(t) &-p_s^{-1}(t)q_s(t)\\
-q_s^*(t)p_s^{-1}(t) &q_s^*(t)p_s^{-1}(t)q_s(t)-r_s(t)
\end{array}\right).
\end{equation}
For each $s\in[0,1]$, let $\Gamma_s(t)$ be a fundamental solution of the
linear Hamiltonian equation

\begin{equation} \label{eq01m1.2}
\dot u=Jb_s(t)u.
\end{equation}

We define the symplectic structure $\omega$ on $\CC^{4m}$ by
$$\omega(x,y) := \lla (-J)\oplus Jx,y\rra,\qquad \forall x,y\in \CC^{4m}.$$

Let $\{W_s\}$ be a continuous family of Lagrangian subspaces of
$\CC^{4m}$. We define
\begin{eqnarray*}
u_s(x)& := &(p_s\frac{d}{dt}x+q_sx,x),\\
D_s& := &\{x\in H^2([0,T],\CC^m)|\bigl((u_s(x))(0),(u_s(x))(T)\bigr)\in W_s\}.
\end{eqnarray*}

\begin{theorem}\label{t:hamiltonian} We define the operators $L_{s,W_s}$ to
be the operators $L_s$ with domain $D_s$. Then they form a continuous
family of self-adjoint Fredholm operators on $L^2([0,T],\CC^m)$. Moreover
we have $$\SF\{L_{s,W_s}\}=\Mas\{{\GGG}\bigl(\Gamma_s(T)\bigr),W_s\}.$$
\end{theorem}

\begin{rem} When all matrices $p_s(t)$ are positive definite, then the
operators $L_s$ are essentially positive. In this case, however,
the spectral flow $\SF\{L_{s,W_s}\}$ is in
general not the same as the difference between the Morse indices
of $L_{0,W_0}$ and $L_{1,W_1}$.

\end{rem}

\subsection{Plan of the Paper}
To prove the four theorems we develop a broad functional analytic frame
work. In Section \ref{s:weak} we develop the symplectic linear algebra and
functional analysis needed, based on a rigorous treatment of (weak)
symplectic vector spaces, Banach spaces, and Hilbert spaces.
We show how to treat continuously varying weak symplectic structures
and define the Maslov index.

In Section \ref{s:sa-extensions} we treat continuous families of closed symmetric operators
in fixed Hilbert space with varying maximal domains but fixed minimal domain
and fixed suitable intermediate domain.
We show that the traces at the boundary of the solution spaces
yield Lagrangian subspaces in corresponding symplectic Hilbert spaces,
if there exist self-adjoint Fredholm extensions.
These {\it reduced Cauchy data spaces} vary continuously under the assumption of
weak inner Unique Continuation Property (UCP).
We discuss various concepts of UCP and prove
the stability of weak inner UCP.
In Subsection \ref{s:proof} we prove the general spectral flow formula in abstract setting.

In Section \ref{s:geometric} we address the geometric setting, i.e., families of symmetric
generalized Dirac type operators over a compact smooth Riemannian
manifold with boundary conditions, which are given by pseudo-differential projections.
We discuss the self-adjointness and the continuity of the corresponding
operators and show that the assumptions of Subsection \ref{s:proof} are naturally
satisfied in our applications.
This proves our two main results, Theorem \ref{t:gsff1} and Theorem \ref{t:split}.
As corollaries we obtain two applications to Hamiltonian dynamics in Subsection \ref{ss:hamiltonian}.

In Appendix \ref{s:appendix1} we give a rigorous definition of the spectral flow in the
unbounded, non self-adjoint case.

\medskip
{\bf Acknowledgement}. We would like to thank Prof. Dr. M. Lesch (K{\"o}ln)
for inspiring discussions about this subject.

\bigskip

\addtocontents{toc}{\medskip\noi}
\section{Weak Symplectic Functional Analysis}\label{s:weak}

\subsection{Basic Symplectic Functional Analysis}\label{ss:basic_symplectic}
We fix our notation. To keep track of the required assumptions we shall
not always assume that the underlying space is a Hilbert space but permit
Banach spaces and -- for some concepts -- even just vector spaces.

\begin{definition}\label{d:symplectic_space}
Let $H$ be a complex Banach space. A mapping
\[
  \w:H\times H\too \C
\]
is called a (weak) symplectic form on $H$, if it is sesquilinear, bounded,
skew-symmetric, and non-degenerate, i.e.,

\noi (i) $\w(x,y)$ is linear in $x$ and conjugate linear in $y$;

\noi (ii) $|\w(x,y)| \leq C\|x\| \|y\|$ for all $x,y\in H$;

\noi (iii) $\w(y,x)=-\ol{\w(y,x)}$;

\noi (iv) $H^{\w} := \{x\in H \mid \w(x,y)=0\text{ for all $y\in H$}\} =
\{0\}$.

\noi Then we call $(H,\w)$ a {\em (weak) symplectic Banach space}.
\end{definition}

There is a purely algebraic concept, as well.
\begin{definition}\label{d:algebraic_symplectic_space}
Let $H$ be a complex vector space and $\w$ a form which satisfies all the
assumptions of Definition \ref{d:symplectic_space} except (ii). Then we call
$(H,\w)$ a {\em complex symplectic vector space}.
\end{definition}

\begin{definition}\label{d:lagrangian}
Let $(H,\w)$ be a complex symplectic vector space.

\noi (a) The {\em annihilator} of a subspace ${\la}$ of $H$ is defined by
\[
{\la}^{\w} := \{y\in H \mid \w(x,y)=0 \text{ for all $x\in {\la}$}\} .
\]

\noi (b) A subspace ${\la}$ is called {\em isotropic}, {\em co-isotropic},
or {\em Lagrangian} if
\[
{\la} \,\<\, {\la}^{\w}\,,\quad {\la}\,\>\, {\la}^{\w}\,,\quad {\la}\,=\,
{\la}^{\w}\,,
\]
respectively.

\noi (c)
The {\em Lagrangian Grassmannian} $\Ll(H,\w)$ consists of all Lagrangian
subspaces of $(H,\w)$.
\end{definition}

\begin{rem}\label{r:lagrangian}
(a)   By definition, each 1-dimensional subspace in real symplectic space
is isotropic, and there always exists a Lagrangian subspace.  However,
there are complex symplectic Hilbert spaces without any Lagrangian subspace.

\noi (b) If $\dim H$ finite, a subspace ${\la}$ is Lagrangian if and only
if it is isotropic with $\dim {\la} = \12 \dim H$.

\noi (c) In symplectic Banach space, the annihilator ${\la}^{\w}$ is
closed for any subspace ${\la}$. In particular, all Lagrangian subspaces
are closed, and we have for any subspace ${\la}$ the inclusion
\begin{equation}\label{e:double_annih}
{\la}^{\w\w} \> \ol{{\la}}.
\end{equation}

\noi (d) Let $H$ be a vector space and denote its (algebraic) dual space by $H'$.
Then each symplectic form $\w$ induces a uniquely defined injective mapping
$J:H\to H'$ such that
\begin{equation}\label{e:almost_complex}
\w(x,y) = (Jx,y) \text{ for all $x,y\in H$},
\end{equation}
where we set $(Jx,y):=(Jx)(y)$.

If $(H,\w)$ is a symplectic Banach space, then the induced mapping $J$ is
a bounded, injective mapping $J:H\to H^*$ where $H^*$ denotes the
(topological) dual space. If $J$ is also surjective (so, invertible), the
pair $(H,\w)$ is called a {\em strong symplectic Banach space}.
We have taken the distinction between {\em weak} and {\em strong} symplectic
structures from Chernoff and Marsden \cite[Section 1.2, pp. 4-5]{ChMa74}.

If $H$ is a Hilbert space with symplectic form $\w$, then the induced
mapping $J$ is a bounded, skew self-adjoint operator (i.e., $J^* =-J$) on
$H$ with $\ker J=\{0\}$.
\end{rem}

\smallskip

The proof of the following lemma is straightforward and is omitted.

\begin{lemma}\label{l:weak-strong}
Any strong symplectic Hilbert space $(H,\lla\cdot,\cdot\rra,\w)$ (i.e., with invertible
$J$) can be made into a
strong symplectic Hilbert space $(H,\lla\cdot,\cdot\rra',\w)$ with $J'^2=-I$ by smooth deformation
of the inner product of $H$ into
\[
\lla x,y\rra':= \lla \sqrt{J^*J}x,y\rra
\]
without changing $\w$.
\end{lemma}

\begin{rem}\label{r:weak-strong}
(a) In a strong symplectic Hilbert space many calculations become quite
easy. E.g., the inclusion \eqref{e:double_annih} becomes an equality, the
Lagrangian property of a subspace ${\la}$ can be characterized by
$(J\la)^\perp={\la}$, and all Fredholm pairs of Lagrangian subspaces have
vanishing index.

\noi (b) We shall give an important example of a weak symplectic Hilbert space:
Let $A$ be a generalized Dirac type operator
in the sense of Definition \ref{d:dirac} over a smooth compact Riemannian manifold $M$
with boundary $\Si$. As mentioned in the Introduction, we have (we suppress mentioning
the vector bundle)
\[
H^{1/2}(\Si) \simeq H^1(M)/H^1_0(M)
\]
with uniformly equivalent norms. Green's form yields a strong symplectic structure
on $L^2(\Si)$ by
\[
\{x,y\} := -\lla Jx,y\rra_{L^2(\Si)}\,.
\]
Here $J$ denotes the principal symbol of the operator $A$ over the boundary
in inner normal direction. It is invertible by our assumption in Definition \ref{d:dirac}.
For the induced symplectic structure on $H^{1/2}(\Si)$  we define $J'$ by
\[
\{x,y\}  = -\lla J'x,y\rra_{H^{1/2}(\Si)}  \quad\text{ for $x,y\in H^{1/2}(\Si)$}.
\]
Let $B$ be a formally self-adjoint elliptic operator $B$ of first order
on $\Si$. By G{aa}rding's inequality, the $H^{1/2}$ norm is equivalent to the
induced graph norm. This yields $J'=(I+|B|)\ii J$.
Since $B$ is elliptic, it has compact resolvent.
So, $(I+|B|)\ii$ is compact in $L^2(\Si)$; and so is $J'$. Hence $J'$ is
not invertible.

\noi (c) Each weak symplectic Hilbert space $(H,\lla\cdot,\cdot\rra,\w)$
with induced injective skew-self-adjoint $J$ can naturally be embedded in a strong symplectic
Hilbert space $H',\lla\cdot,\cdot\rra',\w ')$ with invertible induced $J'$ by setting
$\lla x,y\rra ':= \lla|J| x,y\rra$ as in Lemma \ref{l:weak-strong} and then completing
the space. This imitates the situation of the embedding of $H^{1/2}(\Si)$ into $L^2(\Si)$\,.
It shows that the weak symplectic Hilbert space $H^{1/2}(\Si)$ with its embedding into $L^2(\Si)$
yields a model for all weak symplectic Hilbert spaces.

\end{rem}

\medskip
A key result in symplectic analysis is the following lemma. The
representation of Lagrangian subspaces as graphs of unitary
mappings from one component $H^+$ to the complementary component
$H^-$ of the underlying symplectic vector space (to be considered
as the induced complex space in classical real symplectic
analysis, see, e.g., Booss-Bavnbek and Furutani \cite[Section
1.1]{BoFu98}) goes back to Leray \cite{Le78}. We give a
simplification for complex vector spaces, first announced in
\cite{Zh01}. Of course, the main ideas were already contained in
the real case.

\begin{lemma}\label{l:strong_symplectic}
Let $(H,\w)$ be a strong symplectic Hilbert space with $J^2=-I$. Then
\begin{enumerate}
\item
the space $H$ splits into the direct sum of mutually orthogonal
closed subspaces
\[
H=\ker (J-iI)\oplus\ker(J+iI),
\]
which are both invariant under $J$;

\item there is a 1-1 correspondence between the space $\Uu^J$ of unitary
operators from $\ker(J-iI)$ to $\ker(J+iI)$ and $\Ll(H,\w)$ under the
mapping $U\to {\la}:= \GGG(U)$ (= graph of $U$);

\item if $U,{U'}\in\Uu^J$ and ${\la}:=\GGG(U)$, $\mu:=\GGG({V})$, then
$({\la},\mu)$ is a Fredholm pair (see Definition \ref{d:fredholm_pair}b)
 if and only if $U-V$, or, equivalently,
$UV\ii-I_{\ker(J+iI)}$ is Fredholm. Moreover, we have a natural
isomorphism
\begin{equation}\label{e:unitary_counting}
\ker(UV\ii-I_{\ker(J+iI)}) \simeq {\la}\cap \mu\,.
\end{equation}
\end{enumerate}

\end{lemma}

The proof of (i) is clear; (ii) will follow from Lemma \ref{l:lagrangian_representation};
and (iii) from Proposition \ref{p:fp_characterization}.

\medskip

The preceding method to characterize Lagrangian subspaces and to
determine the dimension of the intersection of a Fredholm pair of
Lagrangian subspaces provides the basis for defining the Maslov
index in strong symplectic spaces of infinite dimensions (see, in
different formulations and different settings, the quoted
references \cite{BoFu98}, \cite{BoFuOt01}, \cite{FuOt02},
\cite{KiLe00}, and Zhu and Long \cite{ZhLo99}).

Surprisingly, it can be generalized to weak symplectic Banach spaces in
the following way.

\begin{lemma}\label{l:lagrangian_representation}
Let $(H,\w)$ be a symplectic vector space and $H^+,H^-$ subspaces. We
assume that $H=H^+ \oplus H^-$ and that the quadratic form $-i\w$ is
positive definite on $H^+$ and negative definite on $H^-$\,. Moreover, we
assume that
\begin{equation}\label{e:h_transversality}
\w(x,y)=0 \text{ for all $x\in H^+$ and $y\in H^-$}\,.
\end{equation}

\noi (a) Then each isotropic subspace ${\la}$ can be written as the graph
\[
{\la}=\GGG(U)
\]
of a uniquely determined injective operator
\[
U:\dom U \too H^-
\]
with $\dom U \< H^+$\,. Moreover, we have
\begin{equation}\label{e:almost_unitary}
\w(x,y) = -\w(Ux,Uy) \text{ for all $x,y\in\dom U$}.
\end{equation}

\noi (b) If $H$ is a Banach space, then the part spaces are always closed
and the operator $U$, defined by a Lagrangian subspace ${\la}$ is closed
as an operator from $H^+$ to $H^-$ (not necessarily densely defined).

\noi (c) For Lagrangian subspaces in a strong symplectic Banach space, we
have $\dom U=H^+$ and $\ran U =H^-$; i.e., the generating $U$ is bounded
and surjective with bounded inverse.
\end{lemma}

\begin{proof} {\it a}. Let ${\la}\< H$ be isotropic and $v_++v_-, w_++w_- \in {\la}$
with $v_\pm,w_\pm\in H^\pm$\,. By the isotropic property of ${\la}$ and
our assumption about the splitting $H=H^+\oplus H^-$ we have
\begin{equation}\label{e:symplectic-splitting}
0=\w(v_++v_-,w_++w_-)=\w(v_+,w_+)+\w(v_-,w_-).
\end{equation}
In particular, we have
\[
\w(v_++v_-,v_++v_-) = \w(v_+,v_+)+\w(v_-,v_-)=0
\]
and so $v_-=0$ if and only if $v_+=0$. So, if the first (resp. the second)
components of two points $v_++v_-,w_++w_-\in {\la}$ coincide, then also
the second (resp. the first) components must coincide.

Now we set
\[
\dom(U) := \{x\in H^+\mid \exists{y\in H^-} \text{ such that $x+y\in
{\la}$}\}.
\]
By the preceding argument, $y$ is uniquely determined, and we can define
$Ux:=y$. By construction, the operator $U$ is an injective linear mapping,
and property \eqref{e:almost_unitary} follows from
\eqref{e:symplectic-splitting}.

\noi {\it b}. One checks easily that $H^\pm=(H^\mp)^\w$\,. Annihilators
are always closed. This proves the first part of ({\it b}). Now let
${\la}$ be a Lagrangian subspace, i.e., ${\la}={\la}^\w$\,. So, ${\la}$ is
closed. It is the graph of $U$. So $U$ is closed.

\noi {\it c}. Now ${\la}$ is a Lagrangian subspace in a strong symplectic
Banach space $H$. First we show that $U$ is densely defined in $H^+$\,.
Indeed, if $\ol{\dom (U)} \neq H^+$\/, there would be a $v\in V$ where $V$
denotes the orthogonal complement of $\dom U$ in $H^+$ with respect to the
inner product on $H^+$ defined by $-i\w$. Clearly $(\dom U)^\w =V+H^-$\,.
So, $V=(\dom U)^\w \cap H^+$\,. Then $v+0\in {\la}^\w\setminus {\la}$.
That contradicts the Lagrangian property of ${\la}$. So, we have $\ol{\dom
(U)} = H^+$\,.

Next we see from \eqref{e:almost_unitary} that $U$ is continuous and so is
$U\ii$\,. Since $U$ has a closed graph, it follows that $\dom U=H^+$ and
$U$ is bounded and has bounded inverse. Applying the same arguments to
$\dom U\ii\<H^-$, relative to the inner product $i\w$ yields $\ran U =\dom
U\ii= H^-$.
\end{proof}

\begin{rem}\label{r:splitting}
(a) Note that the splitting is not unique. Its existence may be proved by
Zorn's Lemma. In our applications, the geometric background provides
natural splittings (see Equation \ref{e:hpm_in_h12}).
For varying splittings see also the discussion below
in Subsection \ref{ss:maslov}.

\noi (b) The symplectic splitting and the corresponding {\em graph} representation
of isotropic and Lagrangian subspaces must be distinguished from the
splitting in complementary Lagrangian subspaces which yields the common
representation of Lagrangian subspaces as {\em images} in the real category (see
below Lemma \ref{l:mas-bofu=mas-zhu}).
\end{rem}

\medskip

\subsection{Fredholm Pairs of Lagrangian Subspaces}
A main feature of symplectic analysis is the study of the {\em Maslov index}.
It is an intersection index between a path of Lagrangian subspaces with the
{\em Maslov cycle}, or, more generally, with another path of Lagrangian subspaces.

Before giving a rigorous definition of the Maslov index in weak symplectic functional
analysis (see below Subsection \ref{ss:maslov}) we fix the terminology and give several
simple criteria for a pair of isotropic subspaces to be Lagrangian.

We recall:

\begin{definition}\label{d:fredholm_pair}
(a) The space of (algebraic) \emph{Fredholm pairs} of
linear subspaces of a vector space $H$ is defined by
\begin{equation}\label{e:fp_alg}
\Ff^{2}_{\operatorname{alg}}(H):=\{({\lambda},{\mu})\mid  \dim\left(
{\lambda}\cap{\mu}\right)  <+\infty\;\;
\text{and} \dim \bigl(H/(\lambda+\mu)\bigr)<+\infty\}
\end{equation}
with
\begin{equation}\label{e:fp_index}
\index(\la,\mu):=\dim(\la\cap\mu) - \dim(H/(\la+\mu)).
\end{equation}

\noi (b) In a Banach space $H$, the space of (topological) \emph{Fredholm pairs} is defined by
\begin{multline}\label{e:fp}
\Ff^{2}(H):=\{({\lambda},{\mu})\in\Ff^2_{\operatorname{alg}}(H)\mid
{\lambda},{\mu}, \tand {\lambda}+{\mu} \subset H \text{ closed}\}.
\end{multline}
\end{definition}

\begin{rem}\label{r:redholm-pairs}
Actually, in Banach space the closedness of $\la+\mu$ follows from its
finite codimension in $H$ in combination with the closedness of $\la,\mu$
(see \cite[Remark A.1]{BoFu99}). So, the set of algebraic Fredholm pairs
of Lagrangian subspaces of a symplectic Banach space $H$ coincides with
the set $\Ff\Ll^2(H)$ of topological Fredholm pairs of Lagrangian
subspaces of $H$.
\end{rem}

We begin with a simple algebraic observation.

\begin{lemma}\label{l:isotropic_sum_to_lagrangian}
Let $(H,\w)$ be a symplectic vector space with transversal subspaces $\la,\mu$\,.
If $\la,\mu$ are isotropic subspaces, then they are Lagrangian subspaces.
\end{lemma}

\begin{proof} {\it a}.
From linear algebra we have
\[
\la^{\w}\cap \mu^{\w} =(\la+\mu)^{\w} =\{0\},
\]
since $\la+\mu=H$. From
\begin{equation}\label{e:isotropics}
\la\<\la^{\w}, \mu\<\mu^{\w}
\end{equation}
we get
\begin{equation}\label{e:new_splitting}
H= \la^{\w} \oplus \mu^{\w}\,.
\end{equation}
To prove $\la^{\w}=\la$ (and similarly for $\mu$),
we consider a $x\in\la^{\w}$\,. It can be written
in the form $x=y+z$ with $y\in\la$ and $z\in \mu$ because of the splitting
$H= \la\oplus\mu$. Applying \eqref{e:isotropics} and the splitting
\eqref{e:new_splitting} we get $y=x$ and so $z=0$, hence $x\in\la$.
\end{proof}

With a little work, the preceding lemma can be generalized from
direct sum decomposition to (algebraic) Fredholm pairs.

\begin{proposition}\label{p:isotropic_fp_to_lagrangian}
Let $(H,\w)$ be a symplectic vector space and $(\la,\mu)
\in \Ff^2_{\operatorname{alg}}(H)$\,.

\noi (a)
If $\la,\mu$ are isotropic subspaces with $\index(\la,\mu)\geq 0$,
then
\[
\index(\la,\mu)= 0 \tand (\la+\mu)^{\w}= \la\,\cap\,\mu\,.
\]

\noi (b) If, moreover, $(\la+\mu)^{\w\w}= \la+\mu$,
then $\la$ and $\mu$ are Lagrangian subspaces of $H$.
\end{proposition}

\begin{proof}
First we show that $\wt{H}:= (\la+\mu)/(\la\cap\mu)$ is a symplectic vector space
with the induced form
\[
\wt{\w}([x+y],[\xi+\eta]):= \w(x+y,\xi+\eta) \text{ for $x,\xi\in\la \tand y,\eta\in \mu$},
\]
where $[x+y]:= x+y+\la\cap\mu$ denotes the class of $x+y$ in $\frac {\la+\mu}{\la\cap\mu}$\,.
Since $\la,\mu$ are isotropic, we have $\w(x+y+z,\xi+\eta+\z)=\w(x+y,\xi+\eta)$ for any $z,\z\in
\la\cap\mu$. So $\wt{\w}$ is well defined and inherits the algebraic properties from $\w$.

To show that $(\wt{H})^{\wt{\w}}=\{0\}$, we observe
\begin{equation}\label{e:symp_alg}
(\la+\mu)^{\w}=\la^{\w}\cap\mu^{\w} \> \la\cap\mu\,.
\end{equation}
We also have
\begin{equation}\label{e:symp_alg_mod}
\dim(\la+\mu)^{\w}=\dim(\la\cap\mu)\,,
\end{equation}
hence
\begin{equation}\label{e:symp_alg_modd}
\frac{(\la+\mu)^{\w}}{\la\cap\mu} = \{0\}.
\end{equation}
To see \eqref{e:symp_alg_mod} we recall from Remark \ref{r:lagrangian}d
that the induced mapping $J:H\to H'$ is injective, so
\[
J|_{(\la+\mu)^{\w}} : (\la+\mu)^{\w} \overset{\simeq}{\too} J\bigl((\la+\mu)^{\w}\bigr).
\]
Since $H/(\la+\mu)$ has finite dimension, we can find
a finite-dimensional subspace $h\< H$ such that $H=h\oplus(\la+\mu)$. So,
any $f\in J\bigl((\la+\mu)^{\w}\bigr)\<H'$ is defined by its values $f(u+v)$ for
$u\in h$ and $v\in \la+\mu$. Clearly, $f(v)=\w(x,v)=0$ for all $v\in \la+\mu$, if $f=Jx$
with $x\in (\la+\mu)^{\w}$\,. So, the functional $f$ is completely determined
by its values on $h$. This implies
\[
\dim (\la+\mu)^{\w} = \dim J\bigl((\la+\mu)^{\w}\bigr) \leq \dim h = \dim
\frac H{\la+\mu} \leq \dim(\la\cap\mu).
\]
Here the last inequality is just the non-negativity
of the Fredholm index as defined in \eqref{e:fp_index}.
This proves \eqref{e:symp_alg_mod}.

\smallskip

\noi {\it b}.  Combining \eqref{e:symp_alg_mod} with \eqref{e:symp_alg} yields
\eqref{e:symp_alg_modd}. Moreover, one checks that
\begin{equation}
\Bigl(\frac{\la+\mu}{\la\cap\mu}\Bigr)^{\wt{\w}}
= \frac{(\la+\mu)^{\w}}{\la\cap\mu}\,.
\end{equation}
With \eqref{e:symp_alg_modd} that proves that $\frac{\la+\mu}{\la\cap\mu}$ is a true symplectic
vector space for the induced form $\wt{\w}$ which is spanned by the transversal isotropic
subspaces
\[
\frac{\la+\mu}{\la\cap\mu} = \frac{\la}{\la\cap\mu} \oplus \frac{\mu}{\la\cap\mu}\,.
\]
By Lemma \ref{l:isotropic_sum_to_lagrangian}, the spaces
$\frac{\la}{\la\cap\mu}\,, \frac{\mu}{\la\cap\mu}$ are Lagrangian subspaces.

Clearly $\la\<\la^{\w}\<\la^{\w}\cap(\la+\mu)$. Now consider $x\in\la$ and $y\in\mu$
with $x+y\in\la^{\w}$\,. Then
\[
[x+y] \in \Bigl(\frac{\la}{\la\cap\mu}\Bigr)^{\wt{\w}} = \frac{\la}{\la\cap\mu}
\]
by the Lagrangian property of $\frac{\la}{\la\cap\mu}$\,. It follows that
$x+y\in\la$, hence
\begin{equation}\label{e:almost_lagrangian}
\la^{\w}\cap(\la+\mu) = \la \text{ and similarly } \mu^{\w}\cap(\la+\mu) = \mu\,.
\end{equation}
Combined with the consequence
\[
\la^{\w} \< (\la\cap\mu)^{\w} = \bigl((\la+\mu)^{\w}\bigr)^{\w} = \la +\mu
\]
of our special assumption, the inclusion $\la\>\la^{\w}$ follows and so
the Lagrangian property of $\la$ (and similarly of $\mu$).
\end{proof}

\begin{rem}
(a) The preceding proposition shows that Fredholm pairs of Lagrangian subspaces in
symplectic Hilbert space can not have positive index. In contrast to the strong case, one
may expect that we have pairs with negative index in weak symplectic Hilbert space.
By now, however, this is an open problem.

\noi (b) The delicacy of Lagrangian analysis in weak symplectic Hilbert space may
also be illuminated by addressing the orthogonal projection onto a Lagrangian subspace.
In strong symplectic Hilbert space with unitary $J$,
the range of an orthogonal projection is Lagrangian if and
only if the projections $P$ and $I-P$ are conjugated by the $J$ operator in the well-known way
\[
I-P=JPJ^*\,.
\]
In weak symplectic analysis, $J$ maps the Range $\range P$ onto a
dense subset of $\ker P$, but there the argument stops.

\noi (c) There is another potential difference between the weak and the strong case:
consider the space $\Ff\Ll_{\la}(H)$ of all Lagrangian subspaces which form a
Fredholm pair with a given Lagrangian subspace $\la$\,. Its topology is presently unknown in
the weak case, whereas we have
\[
\pi_1\bigl(\Ff\Ll_{\la}(H)\bigr)\cong \Z
\]
in strong symplectic Hilbert space $H$ (see \cite[Corollary 4.3]{BoFu99}).
\end{rem}

\medskip

In our applications we typically will have chains of symplectic spaces
like the Sobolev spaces $H:=H^{1/2}(\Si;E|_{\Si})$,
$L:=L^2(\Si;E|_{\Si})$, or the distribution space $\bbb_s:=\dom(A^*_s)/\dom(A_s)$
with symplectic forms all defined by Green's form. We have $H\<L$ dense and
$H\< \bbb_s$ dense, whereas the relation between $L$ and $\bbb_s$ is more difficult
(see the {\it Criss-cross Reduction} of \cite{BoFuOt01} and our Remark \ref{r:op-continuity} below).

Here, a few simple observations can be made how to {\em lift} Lagrangian
subspaces (see also below Lemma \ref{l:maslov-embedding} for parameter
dependence of the lifting).

\begin{lemma}\label{l:embedding}
Let $H$ be a Banach space, $(L,\w_L)$ a (weak or strong) symplectic Banach space, and
$i:H\to L$ a linear continuous injective mapping with dense image. Then we have

\noi (a)
$(H,\w_H)$ is a (weak) symplectic Banach  space where the form
$\w_H:H\times H\to \CC$ is defined by
\begin{equation}\label{e:induced-form}
\w_H(x,y):= \w_L(i(x),i(y))\quad\text{ for all $x,y\in H$}.
\end{equation}

\noi (b) Let $\la\< H$ with $\ol{i(\la)}$
a Lagrangian subspace in $L$ and $\ol{i(\la)} \cap i(H)=i(\la)$. Then $\la$
is a Lagrangian subspace of $H$.
In particular, any subspace $\la\< H$ with $i(\la)$
Lagrangian in $L$ is a Lagrangian subspace of $H$.

\noi (c) Now assume that  $\la\< H$ isotropic in $H$.
Then, clearly, $\ol{i(\la)}$ is an isotropic subspace of $L$.
\end{lemma}

\begin{proof} {\it a}. From the continuity of the inclusion $i$ we obtain that
$\w_L$ bounded implies $\w_H$ bounded. Moreover, we find for the annihilator
\[
H^{\w_H} = i\ii\Bigl((i(H))^{\w_L}\bigr)
\]
and, by the density of $i(H)$ in $L$ and the symplectic property of $(L,\w_L)$
\[
(i(H))^{\w_L} = (\ol{i(H)})^{\w_L} = L^{\w_L}=\{0\}.
\]

\noi {\it b}.
By definition, the image of the $H$-annihilator of $\la$
\[
i(\la^{\w_H}) = \{i(x) \mid x\in H \tand \w_H(x,y)=0
\text{ for all $y\in\la$}\}
\]
is contained in the $L$-annihilator of the image $i(\la)$
\[
\bigl(i(\la)\bigr)^{\w_L} = \{z\mid z\in L \tand \w_L(z,i(y))=0 \text{ for all $y\in\la$}\}.
\]
So,
\[
i(\la^{\w_H}) \< \bigl(i(\la)\bigr)^{\w_L}
= \bigl(\ol{i(\la)}\bigr)^{\w_L}=\ol{i(\la)}
\]
by the Lagrangian property of $\ol{i(\la)}$ in $L$. Hence
$\la^{\w_H}\<\la$.

The opposite inclusion is trivial by \eqref{e:induced-form}.
\end{proof}

\smallskip

We close this subsection with the following characterization of Fredholm pairs.

\begin{proposition}\label{p:fp_characterization}
Let $(H,\w)$ be a symplectic Banach space with a direct sum decomposition
$H^+\,,H^-$ such that $-i\w$ is positive definite
on $H^+$, negative definite on $H^-$, and vanishing on $H^+\times  H^-$\,.
Let $\la,\mu$ be isotropic subspaces. Let $U,V$ denote the generating operators
for $\la,\mu$ (in the sense of Lemma \ref{l:lagrangian_representation}).
We assume that $V$ is bounded and invertible. Then

\noi (a) The space $\mu$ is a Lagrangian subspace of $H$.

\noi (b) Moreover,
\[
(\la,\mu)\in \Ff^2(H) \iff UV\ii-I_{H^-} \text{ Fredholm operator}.
\]

\noi (c) In this case,
\[
\index(\la,\mu)=\index(UV\ii -I_{H^-}).
\]
\end{proposition}

\begin{proof} {\it a}. Let $\mu=\GGG(V)$ be an isotropic subspace of $H$ with
$V:H^+\to H^-$ bounded and invertible. Then also $\mu':=\GGG(-V)$ is
isotropic. We show that $\mu,\mu'$ are transversal in $H$. Then by Lemma
\ref{l:isotropic_sum_to_lagrangian}, $\mu$ (and $\mu'$) are Lagrangian
subspaces. First, from the uniqueness of defining $V$ (see, e.g., the
proof of Lemma \ref{l:lagrangian_representation}a), we have
$\mu\cap\mu'=\{0\}$.

Next, let $x+y$, or, more suggestively, $\begin{pmatrix} x\\ y\end{pmatrix}$ denote any
arbitrary point in $H$ with $x\in H^+$ and $y\in H^-$\,. We set
\[
z := \frac {x+V\ii y}2 \ \tand\ w:= \frac{x-V\ii y}2\,.
\]
Then $z+w=x$ and $z-w=V\ii y$, so
\[
\begin{pmatrix} x\\ y\end{pmatrix}
= \begin{pmatrix} z\\ Vz \end{pmatrix}
+ \begin{pmatrix} w\\ -Vw\end{pmatrix}\,.
\]
This proves $H=\mu\oplus\mu'$\,.

\noi {\it b and c}. Let $\la=\GGG(U)$ and $\mu=\GGG(V)$ with $V$ bounded
and invertible. Let $P_+$, respectively $P_-$\,, denote the projection of
$H=H^+\oplus H^-$ onto the first, respectively, the second factor. Then $$
\la\cap\mu =\Bigl\{\begin{pmatrix}x \\ Vx\end{pmatrix}\mid x\in H^+ \tand
Ux=Vx\Bigr\}. $$ So, $P_2$ induces an algebraic and topological
isomorphism between $\la\cap\mu$ and $\ker \bigl(UV\ii - I_{H^-}\bigr)$.

Now we determine
\begin{align*}
\la+\mu & =\Bigl\{\begin{pmatrix}x \\ Ux\end{pmatrix} + \begin{pmatrix}y
\\ Vy\end{pmatrix} \mid x,y\in H^+\Bigr\}\\ &=\Bigl\{\begin{pmatrix}x' \\
Vx'\end{pmatrix} + \begin{pmatrix}0 \\ z\end{pmatrix} \mid x'\in H^+\tand
z\in \ran(UV\ii-I_{H^-})\Bigr\}\\ &= \mu \oplus \ran(UV\ii-I_{H^-}).
\end{align*}
The last direct sum sign comes from the invertibility of $V$ which induces
$\mu\cap H^-=\{0\}$ and, similarly, $\mu+H^-= H$, and so finally the
direct sum decomposition $H=\mu\oplus H^-$ with projections $\Pi_{\mu}$
and $\Pi_-$ onto the components. So, $\Pi_-$ yields an algebraic and
topological isomorphism of $\la+\mu$ onto $\ran(UV\ii-I_{H^-})$. In
particular, we have $\la+\mu$ closed in $H$ if and only if
$\ran(UV\ii-I_{H^-})$ closed in $H^-$\, and
\[
H/(\la+\mu) \simeq H^-/\ran(UV\ii-I_{H^-})
\]
with coincidence of the codimensions.
\end{proof}


\subsection{Spectral Flow for Curves of ``Unitary" Operators}
\label{ss:isometric}
Let $H$ be a complex Banach space. First let
us introduce some notation for various spaces of operators in $H$:
\[%
\begin{array}
[c]{rl}%
\mathcal{C}(H):= & \text{closed operators on $H$},\\
\mathcal{B}(H):= & \text{bounded linear operators $H\to H$},\\
\mathcal{K}(H):= & \text{compact linear operators $H\to H$},\\
\mathcal{F}(H):= & \text{bounded Fredholm operators $H\to H$},\\
\mathcal{CF}(H):= & \text{closed Fredholm operators on $H$}.
\end{array}
\]

If no confusion is possible we will omit ``$(H)$'' and write $\mathcal{C},$
$\mathcal{B},\mathcal{K},$ etc.. By $\mathcal{C}^{\operatorname{sa}%
},\mathcal{B}^{\operatorname{sa}}$ etc., we denote the set of self--adjoint
elements in $\mathcal{C}, $ $\mathcal{B},$ etc..

We assume that $H$ is a pre-Hilbert space, i.e.,
we are given a fixed inner product (i.e., a sesquilinear,
self-adjoint positive definite form) $h:H\times H\to\C$ which is bounded
\[
|h(x,y)|\leq c\norm{x}\norm{y}\text{ for all $x,y\in H$}.
\]

\begin{definition}\label{d:h_unitary}
An operator $A\in\Cc(H)$ will be called {\em unitary with respect to} $h$, if
\[
h(Ax,Ay)=h(x,y) \text{ for all $x,y\in \dom A$}.
\]
\end{definition}

\begin{rem}\label{r:h_unitary}
(a) Note that $h$ induces a uniformly smaller norm on $H$ which makes $H$ into a
Hilbert space if and only if $H$ becomes complete for this $h$-induced norm.

\noi (b) The concept of $h$-unitary extends trivially to closed operators with
dense domain in one Banach space, equipped with an inner product, and range in
a second Banach space, possibly with a different inner product. Exactly in this sense,
for any Lagrangian subspace the generating operator $U\in\Cc(H^+,H^-)$ (established in Lemma
\ref{l:lagrangian_representation}) is $h$-unitary with $h(x,y)=\mp i\w(x,y)$ on $H^\pm$\,.
\end{rem}

Like for unitary operators in Hilbert space, the following lemma shows that
a unitary operator with respect to $h$ has no
eigenvalues outside the unit circle.

\begin{lemma}\label{l:spectrum_unitary}
Let $A\in\Cc(H)$ be unitary with respect to $h$ and $\la\in\C$, $|\la|\neq 1$. Then
$\ker(A-\la I)=\{0\}$.
\end{lemma}

\begin{proof}
Let $x\in\ker(A-\la I)$, so $Ax=\la x$ and
\[
h(x,x)=h(Ax,Ax)=|\la|^2h(x,x).
\]
Since $|\la|\neq 1$, we get $h(x,x)=0$ and so $x=0$ by $h$ positive definite.
\end{proof}

For a certain subclass of unitary operators with respect to $h$ we show
that they have discrete spectrum close to 1. Consequently, they
are admissible with respect to the positive half-line $\ell$
(in the sense of Definition \ref{d:admissible} of our Appendix) and so
permit the definition of spectral flow through $\ell$ for
continuous families (same Appendix).

\begin{proposition}\label{p:unitary_admissible}
(a) Let $H$ be a Banach space with bounded inner product $h$. Let $A\in\Cc(H)$
be unitary with respect to $h$. We assume
$A-I\in\Cc\Ff(H)$ of index 0. Then there is a neighbourhood $N\<\C$ of $1$
such that the geometric and the algebraic multiplicity of $\s(A)\cap N$
are finite and coincide, and
\[
\s(A)\cap N \< \{1\}.
\]

\noi (b) Let $\{h_s\}$ be a continuous family of inner products for $H$.
Let $A_s\in \Cc(H)$ be unitary with respect to $h_s$. We assume that the
family $\{A_s\}$ is continuous. We denote $h_0=: h$ and $A_0=:A$ and
choose $N$ like in (a). We assume that $N$ is compact with smooth
boundary. Then for $s\ll 1$ the spectrum part $\s(A_s)\cap N$ has finite
algebraic multiplicity and we have
\[
\s(A_s)\cap N \< S^1\,.
\]
\end{proposition}

\begin{proof} {\it a}. We exploit the fact that
the space of closed (generally unbounded) Fredholm operators is
open in $\Cc(H)$ in the graph (= gap) norm and the index is
locally constant (for a full proof see H. O. Cordes and J. P.
Labrousse \cite[Theorem 5.1]{CoLa63}),

Since $\ker(A-I)$ is finite-dimensional, we have an $h$-orthogonal splitting
\[
H=\ker(A-I)\oplus H_1
\]
with closed $H_1$\,. (Take $H_1:=\Pi(H)$ with $\Pi(x):=x-\sum_{j=1}^n h(x,e_j)e_j$, where
$\{e_j\}$ is an $H$-orthonormal basis of $\ker(A-I)$). We notice that $\ker(A-I)\<\dom A$, so
\begin{equation}\label{e:h_orthogonal_domain}
\dom A = \ker(A-I) \oplus (\dom A\cap H_1) .
\end{equation}
Consider the operator $A|_{\dom A\cap H_1}$\,. Let $y\in \dom A\cap H_1$\, and $x\in\ker(A-I)$.
Then
\[
h(x,Ay)=h(Ax,Ay)=h(x,y)=0
\]
by \eqref{e:h_orthogonal_domain}. So, the range $\ran(A|_{\dom A\cap H_1})$ is $h$-orthogonal to
$\ker(A-I)$ and, hence, contained in $H_1$\,. So, the operator $A$ can be split into the form
\begin{equation}\label{e:a_restricted_diagonal}
A=I|_{\ker(A-I)}\oplus A|_{\dom A\cap H_1} = \begin{pmatrix} I_0&0\\0&A_1\end{pmatrix}\,,
\end{equation}
where $I_0$ denotes the identity operator on $\ker(A-I)$ and
$A_1$ denotes the restriction $A_1:=A|_{\dom A\cap H_1} :\dom A\cap H_1\to H_1$\,.

We observe that $A-I$ is closed as bounded perturbation of the
closed operator $A$; it follows that the component $A_1$ and the
operator $A_1-I_1$ are closed in $H_1$\,. Here $I_1$ denotes the
identity operator on $\dom A\cap H_1$\,. By construction,
$\ker(A_1-I_1)=\{0\}$. By our assumption, we have $\index(A_1-I_1)
=\index(A-I)=0$, so $A_1-I_1$ surjective. By the Closed Graph
Theorem, it follows that $(A_1-I_1)\ii$ is bounded and so
$A_1-I_1$ has bounded inverse. The same is true for all operators
in a small neighborhood. Hence, $A_1$ has no spectrum near 1. From
the decomposition \eqref{e:a_restricted_diagonal} we get
$\s(A)=\s(I_0)\cup\s(A_1)$ with $\s(I_0)=\{1\}$. So, if
$1\in\s(A)$ it is an isolated point of $\s(A)$ of finite
multiplicity.

{\it b}. From our assumption it follows that $\s(A)\cap \partial N = \emptyset$, and,
actually, $\s(A_s)\cap \partial N = \emptyset$ for $s$ sufficiently small.
Then
\[
P_N(A_s) := -\frac 1{2\pi i}\int_{\partial N} (A-\la I)\ii d\la
\]
is a continuous family of projections. From T. Kato \cite[Lemma
I.4.10]{Ka76} we get
\[
\dim\ran P_N(A_s) = \dim\ran P_N(A) < +\infty \tand
P_N(A_s)A_s\< A_sP_N(A_s),
\]
and from \cite[Lemma III.6.17]{Ka76} we get
\[
\s(A_s)\cap N=\s(P_N(A_s) A_s P_N(A_s)).
\]
Since all $P_N(A_s) A_s P_N(A_s)$ are unitary with respect to $h_s|_{\ran
P_N(A_s)}$, it follows $\s(P_N(A_s) A_s P_N(A_s))\< S^1$\,.
\end{proof}

Thus, it follows that any $h$-unitary operator $A$ with $A-I$ Fredholm of index 0
has the same spectral properties near $|\la|=1$ as unitary operators in Hilbert space
with the additional property that
1 is an isolated point of the spectrum of finite multiplicity. This now permits us to define the
Maslov index in weak symplectic analysis.


\subsection{Maslov Index in Weak Symplectic Analysis}
\label{ss:maslov}
Our data for defining the Maslov index are a {\em continuous} family
\[
\{(H,\w_s, H_s^+,H_s^-)\}
\]
 of weak symplectic Banach spaces with {\it continuous}
splitting and a {\em continuous} family
$\{(\la_s,\mu_s)\}$ of Fredholm pairs of Lagrangian subspaces of $\{(H,\w_s)\}$
of index 0. Our main task is defining the involved ``continuity".

\begin{definition}\label{d:continuous_banach}
Let $H$ be a fixed complex Banach space and $\{\w_s\}$ a family of weak symplectic forms
for $H$. Let $H_s^+$, $H_s^-$ be subspaces of $H$ such that
$H=H_s^+\oplus H_s^-$ with $\mp i\w_s|_{H_s^\pm}$ positive
definite and $\w_s(x,y)=0$ for all $x\in H_s^+$ and $y\in H_s^-$, $s\in [0,1]$.

\noi (a) The family $\{(H,\w_s,H_s^+,H_s^-)\}$ will be called {\em continuous} if
the induced injective mappings $J_s:H\to H^*$ are continuous as bounded operators, and
the families $\{H_s^\pm\}$ are continuous as closed subspaces of $H$ in the gap topology.
Equivalently, we may demand that the family $\{P_s\}$ of projections
\[
P_s: x+y \mapsto x \text{ for $x\in H_s^+ \tand y\in H_s^-$}
\]
is continuous.

\smallskip

\noi (b)  Let $\{(H,\w_s, H_s^+,H_s^-)\}$ be a continuous family
of symplectic splittings and $\{(\la_s,\mu_s)\}$ a continuous curve of Fredholm pairs
of Lagrangian subspaces of index 0.
Let $U_s:\dom U_s\to H^-_s$, resp. $V_s:\dom V_s\to H^-_s$ be closed
$h_s$-unitary operators with $\GGG(U_s)=\la_s$ and $\GGG(V_s)=\mu_s$\,.
We define the {\em Maslov index} of the curve $\{\la_s,\mu_s\}$ with respect to $P_s$ by
\begin{equation}\label{e:maslov}
\Mas\{\la_s,\mu_s;P_s\} := \SF_{\ell}\Bigl\{\begin{pmatrix} 0&U_s\\ V_s\ii&0\end{pmatrix}\Bigr\}
\/,
\end{equation}
where $V\ii$ denotes the algebraic inverse of the closed injective
operator $V$ and $\ell:=(1-\e,1+\e)$ with suitable real $\e>0$ and with
upward co-orientation. The discussion around Lemma \ref{l:boxplus} below
shows that the spectral flow on the right side of \eqref{e:maslov} is
always well defined.
\end{definition}

\begin{rem}\label{r:full_domain}
Let $\{(H,\w_s,H_s^+,H_s^-)\}$ be a continuous family.
A curve $\{\la_s\}$ of Lagrangian subspaces is continuous
(i.e., $\{\la_s=\GGG(U_s\}$ is continuous as a curve of closed subspaces of $H$), if and only
if the family $\{S_{s,s_0}\circ U_s\circ S_{s,s_0}\ii\}$ is continuous as
a family of closed, generally unbounded operators in the space $\ran P_{s_0}$.
Here $U_s$ denotes the generating operator
$U_s:\dom U_s\to H_s^-$ with $\mathfrak{G}(U_s)=\la_s$ (see Lemma \ref{l:lagrangian_representation});
$s_0\in [0,1]$ is chosen arbitrarily to fix the domain of the family; and
\[
S_{s,s_0}:\ran P_s\too \ran P_{s_0}
\]
is a bounded operator with bounded inverse which is defined in the following
way (see also \cite[Section I.4.6, pp. 33-34]{Ka76}):
\[
S_{s,s_0}:= S_{s,s_0}' (I-R)^{-1/2} = (I-R)^{-1/2} S_{s,s_0}'\,,
\]
where
\[
R
:=(P_s-P_{s_0})^2 \ \tand\ S_{s,s_0}' :=P_{s_0}P_s+(I-P_{s_0})(I-P_s).
 \]

\end{rem}

\smallskip

We have the following lemma:

\begin{lemma}\label{l:boxplus}
Let $(H,\w)$ be a weak symplectic Banach space.
Let $\D$ denote the diagonal (i.e., the canonical Lagrangian) in the product symplectic space
$H\boxplus H:=(H,\w)\oplus (H,-\w)$, and $\la,\mu$ are Lagrangian subspaces of $(H,\w)$. Then
\[
(\la,\mu)\in\Ff\Ll^2(H)
\iff (\la\boxplus\mu,\D)\in\Ff\Ll^2\bigl(H\boxplus H\bigr)
\]
and
\[
\index(\la,\mu)=\index(\la\boxplus\mu,\D),
\]
where $\la\boxplus\mu:=\{(x,y)\mid x\in \la,y\in\mu\}$.
\end{lemma}

\begin{proof}
Clearly $(\la\boxplus \mu)\cap\D\simeq \la\cap\mu$,
and $\la\boxplus \mu$, $\D$ are Lagrangian subspaces of $H\boxplus H$.
Since
\[
(\la\boxplus\mu+\D)\cap (\{0\}\boxplus H)\simeq \{0\}\boxplus (\la+\mu),
\]
we have $\la+\mu$ closed, if $\la\boxplus\mu+\D$ is closed.
Re-arranging
\[
\la\boxplus\mu+\D =
\{(x,y)+(\xi,\xi)\mid x\in\la,y\in\mu,\xi\in H\}
=\{(x,y)\mid x-y\in\la+\mu\}
\]
proves the opposite implication. Moreover, we obtain
$\la\boxplus\mu+\D = \D+\D'_{\la+\mu}$ with
$\D'_{\la+\mu}:= \{(x,-x)\mid x\in\la+\mu\}$. So $\la\boxplus\mu +\D$
is closed, if $\la+\mu$ is closed.

Setting, similarly,
$\D' := \{(x,-x)\mid x\in H\}$ yields
\[
\frac{H\boxplus H}{(\la\boxplus \mu)+\D}
= \frac{\D\oplus \D'}{\D\oplus \D'_{\la+\mu}}
\simeq \frac{\D'}{\D'_{\la+\mu}} \simeq \frac{H}{{\la+\mu}}.
\]
This proves our assertion.
\end{proof}

Let $(H,\w)$ be a weak symplectic Banach space with transversal splittings
$H=H_s^+\oplus H_s^-$
and a corresponding  projection $P:H\to H^+$.
Let $(\la,\mu)\in\Ff\Ll^2(H)$. We denote
the generating operators by $U$, respectively $V$\/.
Then we have
\[
\begin{pmatrix} 0&U\\ V\ii&0\end{pmatrix}
=\begin{pmatrix} U&0\\0& V\ii\end{pmatrix}
\begin{pmatrix} 0&I_{H^+}\\ I_{H^-}& 0 \end{pmatrix},
\]
and
\[
\wt\GGG\begin{pmatrix}U&0\\0&V\ii\end{pmatrix}=\la\boxplus\mu,\tand
\wt\GGG\begin{pmatrix}0&I_{H^-}\\I_{H^+}&0\end{pmatrix}=\D,
\]
where $\wt\GGG$ denotes the graph of closed operators from
$\range \Pp$ to $\range (I-\Pp)$ with $\Pp:=P\boxplus (I-P)$.

This leads to the following important result.

\begin{proposition}\label{p:boxplus}
Let $\{(H,\w_s)\}$ be a continuous family of symplectic forms for $H$ in the sense of
Definition \ref{d:continuous_banach}a with transversal splittings $H=H_s^+\oplus H_s^-$
and a corresponding family of projections $\{P_s:H\to H_s^+\}$.
Let $\{(\la_s,\mu_s)\}$ be a continuous curve in $\Ff\Ll^2(H)$. We denote
the generating operators by $U_s$, respectively $V_s$\/.

\smallskip

\noi (a)
If $V_s$ is bounded and has bounded inverse for each $s\in [0,1]$,
then we have
\begin{equation}\label{e:maslov2}
\Mas\{\la_s,\mu_s; P_s\}=\SF_{\ell}\{U_sV_s\ii\},
\end{equation}
where $\ell:=(1-\e,1+\e)$ with suitable real $\e>0$ and with upward co-orientation.

\noi (b) We have
\begin{align}\label{e:maslov1}
\Mas&\{\la_s\boxplus\mu_s,\D;{\Pp}_s\}
=\Mas\{\la_s,\mu_s;P_s\}\\
&=\Mas\{\mu_s,\la_s;I-P_s\}\; \quad\qquad\text {in}\;(H,-\w_s)
\label{e:maslov3}\\
&=\Mas\{\D,\la_s\boxplus\mu_s;I-{\Pp}_s\}\; \quad\text{in}\;(H,-\w_s)\boxplus (H,\w_s),
\label{e:maslov4}
\end{align}
where $\Pp_s:=P_s\boxplus (I-P_s)$\,.
\end{proposition}

\begin{proof}
By our assumption, we have
\[
\dim\ker (z^2I-U_sV_s\ii)=\dim\ker\begin{pmatrix}z I&U_s\\ V_s\ii&z I\end{pmatrix},
\]
for all $z\in\C$. By definition, we have
\[
\Mas\{\la_s,\mu_s; P_s\}
=\SF_{\ell}\Bigl\{\begin{pmatrix} 0&U_s\\ V_s\ii&0\end{pmatrix}\Bigr\}=\SF_{\ell}\{U_sV_s\ii\}.
\]

\noi {\it b}. Let $\wt\GGG$ denote the graph of closed operators from
$\range \Pp_s$ to $\range (I-\Pp_s)$\,. By (a), (b) and c(i) we have
\begin{align*}
\Mas\{\la_s\boxplus\mu_s,\D;{\Pp}_s\}
&= \Mas \Bigl\{\wt\GGG\begin{pmatrix}U_s&0\\0&V_s\ii\end{pmatrix},
\wt\GGG\begin{pmatrix}0&I_{H_s^-}\\I_{H_s^+}&0\end{pmatrix};
\Pp_s\Bigr\} \\
&= \SF_{\ell}\Bigl\{\begin{pmatrix} U_s&0\\0& V_s\ii\end{pmatrix}
\begin{pmatrix} 0&I_{H_s^+}\\ I_{H_s^-}& 0 \end{pmatrix}\Bigr\}\\
&= \SF_{\ell}\Bigl\{\begin{pmatrix} 0&U_s\\ V_s\ii&0\end{pmatrix}\Bigr\}\\
&=\Mas\{\la_s,\mu_s; P_s\}.
\end{align*}
So (\ref{e:maslov1}) is proved.
By the definition of the Maslov index we have (\ref{e:maslov3}).
(\ref{e:maslov4}) follows from (\ref{e:maslov3}) and (\ref{e:maslov1}).
\end{proof}

From the properties of our general spectral flow, as observed at
the end of our Appendix, we get all the basic properties of the
Maslov index (see S. E. Cappell, R. Lee, and E. Y. Miller
\cite[Section 1]{CaLeMi94} for a more comprehensive list).

\begin{proposition}\label{p:maslov-properties}
(a) The Maslov index is {\em invariant under homotopies} of curves of Fredholm pairs
of Lagrangian subspaces with fixed endpoints.
In particular, the Maslov index is invariant under {\em re-parametrization} of
paths.

\noi (b) The Maslov index is additive under {\em catenation},
i.e.
\[
\Mas\bigl\{{\la}_1* {\la}_2, \mu_1 * \mu_2;P_s*Q_s\bigr\}
= \Mas\bigl\{{\la}_1,\mu_1;P_s\bigr\}
+ \Mas\bigl\{{\la}_2,\mu_2;Q_s\bigr\}\,,
\]
where $\{{\la}_i(s)\},\{{\mu}_i(s)\},\, i=1,2$ are continuous paths
with ${\la}_1(1)={\la}_2(0)$, ${\mu}_1(1)={\mu}_2(0)$ and
\[
({\la}_1*{\la}_2)(s) := \begin{cases}
{\la}_1(2s) \quad & 0\leq s\leq \12\\
{\la}_2(2s-1) \quad & \12 < s\leq 1\,,
\end{cases}
\]
and similarly $\mu_1*\mu_2$ and $\{P_s\}*\{Q_s\}$\,.

\noi (c) The Maslov index is {\em natural} under symplectic action: let
$\{(H',\w_s')\}$ be a second family of symplectic Banach spaces and let
\begin{multline*}
L_s\in \operatorname{Sp}(H,\w_s;H',\w_s')\\
:=\{L\in \Bb(H,H')\mid L \text{ invertible and $\w_s'(Lx,Ly)=\w_s(x,y)$}\}
\end{multline*}
such that $\{L_s\}$ is a continuous family as bounded operators. Then, clearly,
$\{H'=L_s(H_s^+)\oplus L_s(H_s^-)\}$ is a continuous family of symplectic splittings
of $\{(H',\w_s')\}$ inducing projections $\{Q_s\}$, and we have
\[
\Mas\{\la_s,\mu_s;P_s\}= \Mas\{L_s\la_s,L_s\mu_s;Q_s\}.
\]

\noi (d) The Maslov index {\em vanishes}, if $\dim \la_s \cap
\mu_s$ constant for all $s\in [0,1]$.

\noi (e) {\em Flipping}. We have
\[
\Mas\{\la_s,\mu_s;P_s\}+ \Mas\{\mu_s,\la_s;P_s\} = \dim \la_0\cap\mu_0
-\dim\la_1\cap\mu_1\/.
\]
\end{proposition}

\medskip

We can not claim that the Maslov index,
$\Mas\{\la_s,\mu_s;P_s\}$ is always independent of
the splitting projection $P_s$ in general Banach space. However, we have the following result.

\begin{proposition}\label{p:mas_independence}
Let $\{(H,\w_s)\}$ be a continuous family of {\em strong} symplectic Banach spaces
(with fixed underlying Banach space $H$) and let $\{H=H_{s,t}^+\oplus H_{s,t}^-\}$ be
two continuous symplectic splittings in the sense of Definition \ref{d:continuous_banach}a
with projections $P_{s,t}:H\to H_{s,t}^+$ for $s\in [0,1]$ and $t=0,1$.
Let $\{(\la_s,\mu_s)\}$ be a continuous curve of Fredholm pairs
of Lagrangian subspaces of $\{(H,\w_s)\}$. Then

\noi (a) \quad $\index (\la_s,\mu_s)=0$ for all $s\in [0,1]$; and

\noi (b) \quad $\Mas\{\la_s,\mu_s;P_{s,0}\} =  \Mas\{\la_s,\mu_s;P_{s,1}\}$ .
\end{proposition}

\begin{note} Commonly, one assumes $J^2=-I$ in strong symplectic analysis
and defines the Maslov index with respect to the induced decomposition.
In view of Lemma \ref{l:weak-strong}, the point of the preceding proposition is
that the Maslov index is independent of the choice of the metric.
\end{note}

\begin{proof} {\it a}. Using $-i\w_s$, we make $(H,\w_s)$ into a symplectic Hilbert
space and deform the metric such that $J_s^2=-I$. Clearly, the dimensions entering into the
definition of the Fredholm index do not change under the deformation. So, we are in the
well-studied standard case.

\noi {\it b}. We recall that our two families of symplectic splitting define
two families of Hilbert structures for $H$ defined by
\begin{multline*}
\lla x,y\rra_{s,t} := -i\w_s(x_{s,t}^+,y_{s,t}^+) + i\w_s(x_{s,t}^-,y_{s,t}^-) \\
\text{for $x=x_{s,t}^+ + x_{s,t}^-, y=y_{s,t}^+ + y_{s,t}^-,
x_{s,t}^+, y_{s,t}^+ \in H_{s,t}^+, x_{s,t}^-, y_{s,t}^- \in H_{s,t}^-$}, t=0,1.
\end{multline*}
For any $t\in [0,1]$ we define
\[
\lla x,y\rra_{s,t} := (1-t)\lla x,y\rra_{s,0} + t\lla x,y\rra_{s,1}\/.
\]
Then all $(H,\lla\cdot,\cdot\rra_{s,t})$ are Hilbert spaces.

Define $J_{s,t}$ by $\w_s(x,y)=\lla J_{s,t}x,y\rra_{s,t}$ and let
$H_{s,t}^\pm$ denote the positive (negative) space of $-iJ_{s,t}$ and $P_{s,t}$
the orthogonal projection of $H$ onto $H_{s,t}^+$\/.

Then the two-parameter family $\{J_{s,t}\}$ is a continuous family of invertible
operators; $\{P_{s,t}\}$ is continuous; and $\{H_{s,t}^+\}$ is continuous.
So $\Mas\{\la_s,\mu_s;P_{s,t}\}$ is well defined. So, by homotopy invariance
and additivity under catenation we obtain
\[
\Mas\{\la_s,\mu_s;P_{s,0}\} = \Mas\{\la_s,\mu_s;P_{s,1}\}.
\]
\end{proof}

\smallskip

For fixed strong symplectic Hilbert space $H$, choosing one single Lagrangian subspace
$\la$ yields a decomposition $H=\la\oplus J\la$\,. This decomposition was used in
\cite[Definition 1.5]{BoFu98} (see also \cite[Theorem 3.1]{BoFuOt01} and
\cite[Proposition 2.14]{FuOt02}) to give the first functional analytic definition
of the Maslov index, though under the somewhat restrictive (and notationally quite demanding)
assumption of {\it real} symplectic structure. Up to the sign,
our Definition \ref{d:continuous_banach}b is a true generalization of that previous definition.
More precisely:

Let $(H,\w)$ be a real symplectic Hilbert space with
\[
\w(x,y)=\lla Jx,y\rra,\ J^2=-I, \ J^t=-J.
\]
Clearly, we obtain a symplectic decomposition $H^+\oplus H^- = H\otimes\C$
with the induced complex strong symplectic form $\w_{\C}$ by
\[
H^\pm:= \{(I \mp\sqm1 J)\z\mid \z\in H\}.
\]

Now we fix one (real) Lagrangian subspace $\la\< H$. Then
there is a real linear isomorphism
$\f:H\cong\la\otimes\C$ defined by $\f(x+Jy)=x+\sqm1 y$ for all $x,y\in\la$.
For $A=X+JY:H\to H$ with $X,Y:H\to H$ real linear and
\begin{equation}\label{e:real-complex-def0}
X(\la)\< \la,\ Y(\la)\< \la,\tand XJ=JX,\ YJ=JY,
\end{equation}
we define
\[
\f_*(A):=\f\circ A\circ \f\ii
=X+\sqm1 Y, \,\ol A_{\la}:=X-JY,\, A^{t_{\la}}:=X^t+JY^t,\,
\]
where $X^t,Y^t$ denotes the real transposed operators.

\begin{lemma}\label{l:mas-bofu=mas-zhu}
Let $(\la,\mu)$ be any pair of Lagrangian subspaces of $H$
(in the real category).
Let $\wt V:H\to H$ with $\wt VJ=J\wt V$ be a real generating operator
for $\mu$ with respect to the orthogonal splitting $H=\la\oplus J\la$,
i.e., $\mu=\wt V(J\la)$ and $\f_*(\wt V)$ is unitary.
Let $U,V:H^+\to H^-$ denote the unitary generating
operators for $\la\otimes\C$ and $\mu\otimes\C$, i.e., we have
\[
\la\otimes \C=\GGG(U)\ \tand \ \mu\otimes \C=\GGG(V).
\]
Then we have $VU\ii=-\ol{S_{\la}(\wt V)}$, where $S_{\la}(\wt V):=\f_*(\wt V)\f_*\bigl(\wt V^{t_{\la}}\bigr)$
is the complex generating
operator for $\mu\otimes \C$ with respect to $\la$,
as defined by J. Leray in \cite[Section I.2.2, Lemma 2.1]{Le78} and elaborated in the preceding
references.
\end{lemma}

\begin{proof}
We firstly give some notations used later.
For $\z=x+Jy\in H$ with $x,y\in\la$ we define
$
\ol{\z}_{\la}:=\f\ii\bigl(\ol{\f(\z)}\bigr)
=x-Jy$.
Moreover, For $A=X+JY:H\to H$ with $X,Y:H\to H$ real linear with
(\ref{e:real-complex-def0}),
we define $\wt S_{\la}(A):=AA^{t_{\la}}$. Then we have
$S_{\la}(A)=\f_*\bigl(\wt S_{\la}(A)\bigr)$.

Now we give explicit descriptions of $U$ and $V$.
It is immediate that $U$ takes the form
\[
\begin{matrix}
U&:&H^+&\too&H^-\\
\ &\ &(I-\sqm1 J)\z&\mapsto&(I+\sqm1 J)\ol\z_{\la}\ .
\end{matrix}
\]
By the definition of $\wt V$, we have
\[
\mu=\wt V(J\la)=\{2\wt VJx+2\sqm1 \wt VJy\mid x,y\in\la\}.
\]
We shall find $V:(I-\sqm1 J)\z\mapsto (I+\sqm1 J)\z_1$ with $\z,\z_1\in H$
such that $\GGG(V)=\mu\otimes\C$, i.e., we shall find $\z_1$ to $\z=x+Jy$
such that
\begin{equation}\label{e:real-complex}
(I-\sqm1 J)\z + (I+\sqm1 J)\z_1 = 2\wt VJx+2\sqm1 \wt VJy \text{ for all
$x,y\in \la$}.
\end{equation}
Comparing real and imaginary part of \eqref{e:real-complex} yields
$\z+\z_1=2\wt VJx$ and $-\sqm1 J(\z-\z_1)=-\sqm1\wt VJy$, so
\[
\z=\wt V(Jx-y) \tand \z_1=\wt V(Jx+y).
\]
From the left equation we obtain
$\ol\z_{\la}=-\ol{\wt V}_{\la}(Jx+y)$.
Since $\f_*(\wt V)$ is unitary, we obtain from the right side
\[
\z_1=\wt V(Jx+y)=-\wt V\ol{\wt V}_{\la}\ii \ol\z_{\la}
=-\wt V\wt V^t\ol\z_{\la}= -\wt S_{\la}(\wt V)\ol\z_{\la}\/.
\]
This gives
\[
\begin{matrix}
V&:&H^+&\too&H^-\\
\ &\ &(I-\sqm1 J)\z&\mapsto&-(I+\sqm1 J)\wt S_{\la}(\wt V)\ol\z_{\la}\/.
\end{matrix}
\]
So for all $z_1:=(I+\sqm1 J)\z_1$ with $\z_1\in H$, we have
\begin{eqnarray*}
VU\ii z_1&=&-(I+\sqm1 J)\wt S_{\la}(\wt V)\z_1\\
&=&-\wt S_{\la}(\wt V)(I+\sqm1 J)\z_1\\
&=&-\wt S_{\la}(\wt V)(I-\sqm1 J)\ol{\f(\z)}\\
&=&-\ol{\f_*(\wt S_{\la}(\wt V))}(I-\sqm1 J)\ol{\f(\z)}\\
&=&-\ol{S_{\la}(\wt V)}(I+\sqm1 J)\z_1\\
&=&-\ol{S_{\la}(\wt V)}z_1.
\end{eqnarray*}
That is, $VU\ii=-\ol{S_{\la}(\wt V)}$.
\end{proof}

\smallskip
With the preceding notation, we recall from \cite[Definition 1.5]{BoFu98} the definition of the Maslov index
\begin{equation}\label{e:maslov-bf}
\Mas_{\operatorname{BF}}\{\mu_s,\la\}:= \SF_{\ell'}\{S_{\la}(\wt V_s)\}
\end{equation}
of a continuous curve $\{\mu_s\}$ of Lagrangian subspaces in real
symplectic Hilbert space $H$ which make Fredholm pairs with one fixed
Lagrangian subspace $\la$. Here $\ell':=(-1-\e,-1+\e)$ with downward orientation.

\begin{corollary}\label{c:mas-bofu=mas-zhu}
\[
\Mas\{\la\otimes\C,\m_s\otimes \C\}=-\Mas_{\operatorname{BF}}\{\mu_s,\la\}.
\]
\end{corollary}

\begin{proof}
Let $\ell,\ell'$ denote small intervals on the real line close to 1, respectively -1
and give $\ell$ the co-orientation from $-i$ to $+i$ and $\ell'$ vice versa. We denote
by $\ell^-$ the interval $\ell$ with reversed co-orientation.
Then by our definition in \ref{e:maslov}, elementary transformations, the preceding lemma,
and the definition recalled in \eqref{e:maslov-bf}:
\begin{align*}
\Mas\{\la\otimes\C,\m_s\otimes \C\}
&= -\SF_{\ell}\{UV_s\ii\}= -\SF_{\ell^-}\{V_sU\ii\} \\
&= -\SF_{\ell^-}\{-\ol{S_{\la}(\wt V_s)}\}
= -\SF_{\ell}\{-S_{\la}(\wt V_s)\}\\
& = -\SF_{\ell'}\{S_{\la}(\wt V_s)\} =-\Mas_{\operatorname{BF}}\{\mu_s,\la\}.
\end{align*}
\end{proof}

We close this section with discussing the invariance of the Maslov index
by embedding in a larger symplectic space, assuming a simple regularity
condition (see also Lemma \ref{l:embedding} above).

\begin{lemma}\label{l:maslov-embedding}
Let $\{(L,\w_s,L^+_s,L^-_s)\}$ be a continuous family of symplectic splittings for a fixed
(complex) Banach space $L$ and
$\{(\la_s,\mu_s)\in \Ff\Ll^2(L,\w_s)\}$ a continuous curve
with $\index(\la_s,\mu_s)=0$ for all $s\in [0,1]$. Let $H$ be a second Banach space with
a linear embedding $H\into L$ (in general neither continuous nor dense). We assume that
\[
\wt{\w}_s:=\w_s|_{H\times H} \tand H_s^\pm:=L_s^\pm\cap H.
\]
yields also a continuous family $\{(H,\wt{\w}_s,H^+_s,H^-_s)\}$ of symplectic splittings.
Moreover, we assume that $\la_s\cap \mu_s\< H$
and $(\la_s\cap H,\mu_s\cap H)\in \Ff\Ll^2(H,\wt{\w}_s)$ of index 0
for all $s$, and that the pairs make also a continuous curve in $H$.
Then we have
\[
\Mas\{\la_s,\mu_s;P_s\} = \Mas\{\la_s\cap H_s,\mu_s\cap H_s;\wt{P}_s\},
\]
where $P_s$ and $\wt P_s$ denote the projections for the decomposition
$L=L_s^+\oplus L_s^-$, respectively $H=H_s^+\oplus H_s^-$\/.
\end{lemma}

The lemma is an immediate consequence of Lemma \ref{l:sf-embedding} of the appendix.


\addtocontents{toc}{\medskip\noi}
\section{Symplectic Analysis of Symmetric Operators}\label{s:sa-extensions}

\subsection{General Assumptions}
In this section we shall introduce the concept of a {\em continuous}
family $\{A_t\}_{t\in[0,1]}$ of closed symmetric operators in a fixed
Hilbert space $X$ (in Section \ref{s:geometric} we shall give an outline
of necessary changes for {\em continuously} varying Hilbert spaces
$X_t$\/), and {\em continuously} varying self-adjoint Fredholm extensions
$A_{t,D_t}$ with {\em continuously} varying domains $\{D_t\< X_t\}_{t\in
[0,1]}$\,. Our main effort will be explaining the suitable concepts of
{\em continuity} and to specify the minimal assumptions required for the
general spectral flow formula.

The most general setting for that will be a fixed (or canonically
identified) minimal domain $\dom A_t=:D_m$; a fixed (or canonically
identified) reduced (intermediate) domain $D_W$; but varying maximal
domains $D_{\mmax,t}:=\dom A_t^*$\,; and varying domains $D_t$ of the
self-adjoint Fredholm extensions.

For a better understanding, one may think of the situation described in the Introduction (see also
below Section \ref{s:geometric}).

Let $X$ be a complex Hilbert space and $A\in\Cc(X)$ a linear, closed, densely defined
operator in $X$. We assume that $A$ is symmetric, i.e., $A^*\>A$ where $A^*$ denotes
the adjoint operator. We denote the domains of $A$ by $D_m$ (the {\it minimal} domain)
and of $A^*$ by $D_{\mmax}$ (the {\it maximal} domain). We recall from \cite{BoFu98}:

\begin{enumerate}

\item
The space $D_{\mmax}$ is a Hilbert space with the graph inner product
\begin{equation}\label{e:graph_inner_product1}
\lla x,y\rra_{\GGG} := \lla x,y\rra_X + \lla A^*x,A^*y\rra_X
\quad\text{ for $x,y\in D_{\mmax}$}\,.
\end{equation}

\item The space $D_m$ is a closed subspace in the graph norm and the quotient
space $D_{\mmax}/D_m$ is a strong symplectic Hilbert space with the
(bounded) symplectic form induced by Green's form
\begin{equation}\label{e:symplectic_green1}
\{x+D_m,y+D_m\} := \lla A^*x,y\rra_X - \lla x,A^*y\rra_X
\quad\text{ for $x,y\in D_{\mmax}$}\,.
\end{equation}

\item If $A$ admits a self-adjoint Fredholm extension $A^*|_D$ with
domain $D\<X$, then the {\it natural} Cauchy data space $(\ker A^*+D_m)/D_m$
is a Lagrangian subspace of $D_{\mmax}/D_m$\,.

\item Moreover, self-adjoint Fredholm extensions are characterized by the
property of the domain $D$ that $(D+D_m)/D_m$ is a Lagrangian subspace of
$D_{\mmax}/D_m$ and forms a Fredholm pair with $(\ker A^*+D_m)/D_m$\,.

\item We denote the natural projection by
\[
\g:D_{\mmax}\too D_{\mmax}/D_m\/.
\]
\end{enumerate}

\smallskip

\subsection{Lagrangian Property of Reduced Cauchy Data
Spaces}\label{ss_symmetric}
In our applications, we consider families of self-adjoint Fredholm operators
with varying domain and varying maximal domain. To us, there is no natural way to identify
the different symplectic spaces and to define continuity of Lagrangian subspaces
and continuity of symplectic forms in
these varying symplectic spaces. However, in most applications the minimal domain is
fixed and also an intermediate (reduced) Hilbert space $D_W$\,, typically the first Sobolev space. We
shall show that meaningful modifications of the preceding statements can be obtained when we
replace $D_{\mmax}$ by this reduced (intermediate) space $D_W$ under the following assumptions.

\begin{asss}\label{a:reduced-space}
(a) Our data are now four Hilbert spaces with continuous inclusions
\[
D_m\into D_W \into D_{\mmax} \into X,
\]
where the Hilbert space structure is given on $D_{\mmax}$ and $D_m$ by the
graph inner product of a fixed closed densely defined symmetric
operator $A\in\Cc(X)$ with $\dom A=D_m$\,.

\noi (b) We assume that on $D_m$ the graph norm and the norm induced by the Hilbert
space $D_W$ are equivalent.

\noi (c) We assume that the set $D_W$ is dense in $D_{\mmax}$ in
the graph norm.

\noi (d) Finally, we assume that there exists
a self-adjoint Fredholm extension $A_D$ of $A$ with domain $D_m\< D \< D_W$\,.
\end{asss}

Assumption \ref{a:reduced-space}a implies
\begin{align}\label{e:am_bounded}
\norm{x}^2_{D_{\mmax}} :&= \norm{x}^2_{\GGG} = \norm{x}^2_X +
\norm{Ax}^2_X\\
&\leq c\norm{x}^2_{D_W} \text{ for all $x\in D_W$}\,.
\end{align}
In particular, it follows that $A^*|_{D_W} : D_W\to X$ is bounded.

Assumption \ref{a:reduced-space}b implies the opposite estimate to \eqref{e:am_bounded},
namely
\begin{equation}\label{e:am_and_gaarding}
\norm{x}^2_{D_W} \leq c\Bigl( \norm{x}^2_X + \norm{Ax}^2_X\Bigr)
= c\norm{x}^2_{\GGG} \text{ for all $x\in D_m$}\,.
\end{equation}

\begin{lemma}\label{l:weak_symplectic_boundary_space}
Under Assumptions \ref{a:reduced-space}a,b,c
the quotient space $D_W/D_m$ is a weak
symplectic Hilbert space with the symplectic form induced by Green's form
on $D_{\mmax}$\,.
\end{lemma}

\begin{proof}
By Assumption \ref{a:reduced-space}b, $D_m$ is closed in $D_W$, so the quotient is a Hilbert space
and we can apply Lemma \ref{l:embedding}a. For easy reading we repeat the argument
of Lemma \ref{l:embedding}a in full length.

By Assumption \ref{a:reduced-space}a the inclusion $D_W\into D_{\mmax}$ is continuous, hence
the restriction of the Green's form on $D_W$ is also bounded.

It remains to show that the annihilator $(D_W/D_m)^0$  of $D_W/D_m$ is $\{0\}$.
We shall denote the class $x+D_m$ in $D_W/D_m$ of any $x\in D_W$ by $[x]$.
So, let $[x]\in (D_W/D_m)^0$\,. Then we have
\begin{equation}\label{e:annihilator_M}
\{x,y\}=0 \quad\text{for all $y\in D_W$}\,.
\end{equation}
From the continuity of the form $\{\cdot,\cdot\}$ on $D_{\mmax}$ and Assumption
\ref{a:reduced-space}c it follows that equation \eqref{e:annihilator_M} is valid for all $y\in
D_{\mmax}$\,. Since $D_{\mmax}/D_m$ is symplectic, that yields $x\in D_m$\,, so
$[x]=0$ in $D_W/D_m$\,.
\end{proof}

The lemma shows that any intermediate space $D_W$ satisfying Assumptions
\ref{a:reduced-space}a,b,c is big enough to permit a meaningful symplectic analysis
on the reduced quotient space $D_W/D_m$\,. The point of this construction is that
the {\em norm} in $D_W/D_m$ does not come from the graph norm in $D_{\mmax}$ but
from the norm of $D_W$\,. Therefore, it can be kept fixed
even when our operator varies. The {\em symplectic structure} of $D_W/D_m$\,, however,
is induced by Green's form and therefore will change with varying operators.

In the following, we denote the
extension $A^*|_D$ of $A$ with domain $D$ by $A_D$\,. We shall write $A_W$
for $A_{D_W}$\,.

In \cite[Proposition 3.5]{BoFu98}, in the spirit of the classical von-Neumann program,
self-adjoint Fredholm extensions were characterized by the property that their
domains, projected down into the strong symplectic space $\bbb(A):=D_{\mmax}/D_m$ of natural
boundary values, make Fredholm pairs of Lagrangian subspaces with the
natural Cauchy data space $(\ker A^*+D_m)/D_m$\/. Immediately, this does not
help for operator families with varying maximal domain. Surprisingly, however,
the arguments generalize to the weak symplectic space $D_W/D_m$\/.

\begin{proposition}\label{p:cauchy}
Under the Assumptions \ref{a:reduced-space}
the quotient space $D/D_m$ and the {\em reduced Cauchy data space}
$(\ker A_W+D_m)/D_m$ form a Fredholm pair of Lagrangian subspaces of the (weak)
symplectic Hilbert space $D_W/D_m$ with index $0$.
Moreover, it follows that $\range A_W=\range A^*$\,.
\end{proposition}

\begin{proof} Clearly, the spaces
\[
\la:=  D/D_m \tand \mu:= (\ker A_W + D_m)/D_m
\]
are isotropic subspaces in the (weak) symplectic Hilbert space $H:= (D_W/D_m,\{\cdot,\cdot\})$\,. So, by
Proposition \ref{p:isotropic_fp_to_lagrangian} it suffices to show that $\la+\mu$ is closed in
$H$ and that $\la,\mu$ form an (algebraic) Fredholm pair of non-negative index.

\smallskip

First, we consider
\begin{equation}\label{e:pi_astar_j}
D_W \ftoo{A^*\circ j} \ran A_W \ftoo{\pi} \ran A_W/\ran A_D\,,
\end{equation}
where $j$ denotes the dense continuous inclusion $D_W \finto{j} D_{\mmax}$ and $\pi$ denotes
the projection onto the quotient space. Note that the range $\ran A_D$ closed in the basic
Hilbert space $X$ by the Fredholm property, and so also closed in the reduced (intermediate)
space $\ran A_W$\,. Then
\[
D+\ker A_W = \{x\in D_W\mid A^*x\in \ran A_D\}=\ker (\pi\circ A^*\circ j)
\]
is closed in $D_W$ since $j, A^*,$ and $\pi$ are continuous. So, also
$\la+\mu = (D+\ker A_W)/D_m$ is closed in $H=D_W/D_m$\,.

\smallskip

To examine $H/(\la+\mu)$ we notice that the reduced (intermediate) space $\ran A_W$
contains the closed space $\ran A_D$
of finite codimension in $X$. So, the space $\ran A_W/\ran A_D$ to the right in
\eqref{e:pi_astar_j} must be of finite dimension. Then
we apply the purely algebraic identities (note that $A_W:D_W\to \ran A_W$
is surjective):
\begin{align}\label{e:quotient_in_proposition}
H/(\la +\mu )&= (D_W/D_m )/ ( (D+\ker A_W)/D_m )\\
&\simeq D_W/\bigl(D+ \ker A_W\bigr)
\simeq \ran A_W/\ran A_D\,.
\end{align}
That proves $\dim H/(\la +\mu )<\infty$.

\smallskip

Next we consider the algebraic identity
\begin{align*}
\la\cap\mu &=D/D_m\cap(\ker A_W+D_m)/D_m \\
&= (\ker A_D+D_m)/D_m  \simeq \ker A_D/(D_m\cap \ker A_D),
\end{align*}
where the numerator and the denominator in the quotient to the right are of finite dimension, so
$\dim (\la\cap\mu)<+\infty$.

\smallskip

To estimate $\index(\la,\mu)$, we re-write, algebraically,
\[
\la\cap\mu \simeq \ker A_D/(D_m\cap \ker A_D) \simeq \ker A_D/ \ker A.
\]
We consider the numerator to the right,
\begin{align*}
\ker A_D&\simeq X/\ran A_D \simeq X/\ran A^*\ \oplus \ran A^*/\ran A_D\\
&\simeq \ker A\ \oplus \ \ran A^*/\ran A_D\,.
\end{align*}
The transformation can be made since $A_D$ is a self-adjoint Fredholm operator with
closed range of finite codimension,
and so $\ran A^*$ with $\ran A_D \< \ran A^*\< X$ is also of finite codimension and closed.
So, we finally obtained an alternative expression for $\dim \la\cap\mu$, namely
\[
\dim \la\cap\mu = \dim \ran A^*/\ran A_D\,.
\]
Comparing with \eqref{e:quotient_in_proposition} we see that
\[
\codim (\la+\mu )= \dim \ran A_W/\ran A_D \leq \dim \la\cap\mu\,.
\]
This proves $\index(\la,\mu)\geq 0$. So, by Proposition \ref{p:isotropic_fp_to_lagrangian},
the pair $(\la,\mu)$ is truly a Fredholm pair of Lagrangian subspaces with
vanishing Fredholm index. So,
\[
\dim \ran A^*/\ran A_D = \dim \ran A_W/\ran A_D \,.
\]
Since $\ran A_W\<\ran A^*$\,, that proves $\ran A_W=\ran A^*$\,.
\end{proof}

By the above proposition we can prove the somewhat surprising
local stability of weak inner UCP. Firstly we recall the definition.

\begin{definition}\label{d:ucp_general}
(a) Let $X$ be a Banach space with subspaces $D_m\< D_W\< X$ and let $A_W$ be an
unbounded operator in $X$ with $\dom A_W=D_W$\/.
We shall say that the operator $A_W$ satisfies the {\em weak inner
Unique Continuation Property (UCP)} with respect to $D_m$ if and only if
$\ker A_W|_{D_m}=\{0\}$.

\noi (b) Let $X$ be a Hilbert space and $A\in \Cc(X)$ with $\dom A=D_m$ and
$A^*\> A$. We shall say that the operator $A$ (or, colloquially,
the operators $A^*$) satisfies the {\em weak inner
Unique Continuation Property (UCP)} if $\ker A=\{0\}$.
\end{definition}

Note that
\[
\ker A = \ker A^*|_{D_m} = \ker A_W|_{D_m} = D_m\cap \ker A^*\,.
\]

It is well known that weak UCP and weak inner UCP can be
established for a great class of Dirac type operators, see the
first author with Wojciechowski \cite[Chapter 8]{BoWo93}, and the
first author with M. Marcolli and B.-L. Wang \cite{BoMaWa02}.
However, it is not valid for all linear elliptic differential
operators of first order as shown by the Pli{\'s} counter-example
\cite{Pl61}. Moreover, one has various quite elementary examples
of linear and non-linear perturbations which {\it invalidate} weak
inner UCP for Dirac operators. Two such examples are listed in
\cite{BoMaWa02}. In the same paper, however, it was shown that
weak UCP is {\it preserved} under certain `small' perturbations of
Dirac type operators. Here we show a more general and more
elementary result, namely the local stability of weak inner UCP.

 \begin{corollary}\label{c:stable-ucp}
Let $X$ be a Hilbert space.
Let $A\in \Cc(X)$ with $\dom A=D_m$\/, and let $D_m\< D\< D_W\< X$ satisfy Assumptions
\ref{a:reduced-space} with respect to $A$. Let
$\{A_{s,W}:D_W\to X\}$ be a continuous curve of bounded operators with $A_{0,W}= A^*|_{D_W}$\,.
If $A_{0,W}$ satisfies weak inner UCP with respect to $D_m$, then all
$A_{s,W}$ are surjective for $s\ll 1$.
If, moreover,
\begin{equation}\label{e:symmetry-ucp}
\lla A_{s,W}x,y\rra = \lla x,A_{s,W}y\rra \quad\text{ for all $x\in D_m$ and $y\in D_W$}\/,
\end{equation}
then,  for $s\ll 1$, all $A_{s,W}$ satisfy weak inner UCP with respect to $D_m$\/.
\end{corollary}

\begin{proof}
By Assumptions \ref{a:reduced-space}, $\ran A^*|_D$ is closed and is of finite codimension. Since
$\ran A^*|_D\<\ran A^*\< X$, the full range $\ran A^*$ is closed.
Since $A_{0,W}= A^*|_{D_W}$ and $A_{0,W}$ satisfies weak inner UCP with respect to $D_m$\/,
$\ran A^*=X$. By Proposition \ref{p:cauchy}, we have $\ran A_{0,W}=\ran A^*$.
Then $A_{0,W}$ is semi-Fredholm. By Theorem IV.5.17 of \cite{Ka76} we have $\ran A_{s,W}=X$
for $s\ll 1$. This proves the first part of the corollary.

Since \eqref{e:symmetry-ucp} means that $\bigl(A_{s,W}|_{D_m}\bigr)^* \> A_{s,W}$,
the second part is an immediate consequence of the first part.
\end{proof}

\smallskip

\subsection{Continuity of Reduced Cauchy Data Spaces}\label{ss:cd-continuity}
In this subsection we generalize (and partly simplify) the proof of the
continuity of Cauchy data spaces given in \cite[Section 3.3]{BoFu98}.
Instead of the natural Cauchy data spaces $\bigl(D_m+\ker A_s^*\bigr)/D_m$
in the natural (distribution) space $D_{\mmax}/D_m$ (see point (iii) of
the introductory remarks to our Subsection \ref{ss_symmetric}) we consider
the reduced Cauchy data spaces $(D_m+\ker A_{s,W})/D_m$ in the reduced
(function) space $D_W/D_m$\,.


\begin{theorem}\label{t:cd_continuity}
Let $\{A_s:D_m\to X\}_{s\in [0,1]}$ be a family of closed symmetric densely defined operators
in $X$. We assume that
\begin{enumerate}
\item each $A_s$ admits a self-adjoint Fredholm extension with domain $D_s$\,;
\item each $A_s^*$ satisfies weak inner UCP relative to $D_m$\,; and
\item the extensions $A_{s,W}=(A_s)^*|_{D_W}:D_W\to X$ form a continuous family
of bounded operators.
\end{enumerate}
Then the reduced Cauchy data spaces $\bigl(D_m+\ker A_{s,W}\bigr)/D_m$ are
continuously varying in $D_W/D_m$\,. Here $D_m\into D_W \into X$ denote
three Hilbert spaces with continuous embeddings satisfying Assumption
\ref{a:reduced-space} relative $A_0$\,.
\end{theorem}

\begin{rem}\label{r:continuity}
(a) As usual, we define the continuous dependence of a family of
subspaces of a Hilbert space on a parameter by the continuity
of the corresponding orthogonal projections.

\noi (b) The assumption that the family $A_{s,W}=A_s^*|_{D_W}:D_W\to X$
forms a continuous family in the space of bounded operators $\Bb(D_W,X)$
is naturally satisfied in our applications, as we shall see below.
However, even in the case of fixed maximal domain $D_W=D_{s,\mmax}$\,, our
assumption is more restrictive than demanding only continuity in the gap
norm, as explained in Remark \ref{r:op-continuity}.

\noi (c) As mentioned before, when we deal with curves of symmetric
generalized operators of Dirac type on a compact manifold $M$ with
boundary $\Si$, we have $D_W=H^1(M)$, $D_m=H^1_0(M)$, and
$D_W/D_m=H^{1/2}(\Si)$. Then the reduced Cauchy data spaces
$\bigl(D_m+\ker A_{s,W}\bigr)/D_m$ can be identified with the ranges $\ran
Q_s\cap H^{1/2}(\Si)$ of the pseudo-differential {\it Calder{\'o}n
projections} $Q_s$ of the operators $A_s$\,.

\end{rem}

\begin{proof}
[Proof of Theorem \ref{t:cd_continuity}] For shorter writing and better
reading, we shall denote $\ker A_s,W$ by $\sss_s$, and the projection of
$D_W$ onto $D_W/D_m$ by $\g$\,. Note that $\sss_s$ is closed in $D_W$\,,
but not necessarily in $D_{s,\mmax}$ nor in $X$.

To prove the continuity, we need only to consider the local
situation at $s=0$. First we show that $\{\sss_s\}_{s\in I}$
is a continuous family of subspaces of $D_W$; then we
show that $\g(\sss_s)$ is a continuous family in $D_W/D_m$.

\smallskip

We consider the bounded operator
\[
\begin{matrix}
F_s: & D_W &\too & X \oplus \sss_0         \\
\ & x& \mapsto & \left(A^*_s(x), P_0x\right)
\end{matrix}\qquad ,
\]
where $P_0:D_W\to \sss_0$ denotes the orthogonal projection
of the Hilbert space $D_W$ onto the closed subspace $\sss_0$\,.
By definition, the family $\{F_s\}$ is a continuous family of
bounded operators.

Clearly, $F_0$ is injective. Moreover, from weak inner UCP we get
$\ran A_0^*=X$ and from Proposition \ref{p:cauchy} $\ran A_0^*=
\ran A_0^*|_{D_W}$\/. So, the operator
$F_0$ is also surjective. This proves that $F_0$ is invertible with
bounded inverse.
Then all operators $F_s$ are invertible for small $s \geq 0$,
since $F_s$ is a continuous family of operators.

Note that
\begin{equation}\label{e:f-range}
F_0(\sss_0)=F(\sss_s)=0\oplus \sss_0\,.
\end{equation}
These linear spaces are just the same.

We define
\[
\f_s:= F_s\ii\circ F_0 :D_W\ \cong\ D_W
\tand
\f_s\ii= F_0\ii\circ F_s :D_W\ \cong\ D_W
\]
for $s$ small.
From \eqref{e:f-range}, we obtain that
\begin{equation}\label{e:phi}
\f_s(\sss_0)=\sss_s\,.
\end{equation}
More explicitly, let $z\in\f_s(\sss_0)$, $z=\f_s(x)=F_s\ii F_0(x)$ for suitable
$x\in \sss_0$\,. Then $F_0(x)=(0,x)$ and $F_s\ii F_0(x)=z\in \sss_s$\,.
To prove the opposite inclusion, let $z\in \sss_s$ and set $x:=P_0(z)$.
Then $F_0(x)=(0,x)$ and $F_s\ii F_0(x)=z$.

From \eqref{e:phi} we get that
\[
\{P_s:=\f_s P_0\f_s\ii:D_W\too \sss_s\}
\]
is a continuous family of projections onto the solution spaces
$\sss_s$. The projections are not necessarily orthogonal, but can
be orthogonalized and remain continuous in $s$ like in
\cite[Lemma 12.8]{BoWo93}. This proves the continuity of the family
$\{\sss_s\}$ in $D_W$\,.

\smallskip

Now we must show that $\{\g(\sss_s)\}$ is a continuous
family in $D_W/D_m$. This is not proved by the formula
$\g(\sss_s)=\g(\f_s(\sss_0))$ alone. We must modify the
endomorphism $\f_s$ of $D_W$ in such a way that it keeps
the subspace $D_m$ invariant.

To do that, we recall that $D_m$\,, which is
closed in the graph norm, is, by Assumption \ref{a:reduced-space}b, closed
in $D_W$\,. So, $D_m+\sss_0$ is closed in
$D_W$. We define a continuous family of mappings by
\[
\begin{matrix}
\psi_s: D_W = &D_m  +  \sss_0 & \ + \ &
(D_m+\sss_0)^\perp & \
\too \ & D_W\\
\ & x+y & + & z&\mapsto&x+\f_s(y)+z
\end{matrix}
\]
with $\psi_0=\id$. Hence all $\psi_s$ are invertible for $s\ll 1$,
and $\psi_s(D_m)=D_m$ for such small $s$. Hence we obtain a
continuous family of mappings $\{\wt\psi_s:D_W/D_m\to D_W/D_m\}$ with
$\wt\psi_s(\g(\sss_0))=\g(\sss_s)$. From that we obtain a
continuous family of projections as above.
\end{proof}

\begin{rem}
From the preceding arguments it also follows that the
Cauchy data spaces form a differentiable family, if $\{A_s^*|_{D_W}\}$
is a differentiable family.
\end{rem}

\smallskip

\subsection{Proof of the General Spectral Flow Formula}\label{s:proof}
We shall make the following assumptions. They all are natural in our
applications, as we shall see later, in Section \ref{s:geometric}.

\begin{asss}\label{a:continuity-symplectic-splitting}
(a) We are given a family $\{A_s:D_m\to X\}_{s\in [0,1]}$ of closed symmetric densely
defined operators satisfying weak inner UCP
and continuity as a family of bounded operators $\{A_s^*:D_W\to X\}$
with four Hilbert spaces with continuous inclusions
\[
D_m\into D_W \into D_{\mmax,s} \into X.
\]
Here the Hilbert space structures on $D_{\mmax,s}:=\dom(A_s^*)$ and $D_m$ is a closed subspace of
$D_W$.
We assume that on $D_m$ the graph norm and the norm induced by the Hilbert
space $D_W$ are equivalent and that the set $D_W$ is dense in each $D_{\mmax,s}$
with the respective graph norm.

\smallskip

\noi (b) We assume that the Hilbert space $D_W/D_m$ is embedded
in a fixed Hilbert space $L$ with a continuous family of symplectic
splittings
$$(L,\w_s,L_s^+,L_s^-).$$
We assume that
\[
\{x,y\}_s=\w_s(x,y)\quad\text{ for $x,y\in D_W/D_m$}\/,
\]
where $\{\cdot,\cdot\}_s$ is induced by the Green's form of $A_s^*$
(see (\ref{e:symplectic_green1})).
We set $H_s^{\pm}:=L_s^{\pm}\cap D_W/D_m$ and assume that
$(D_W/D_m,\{\cdot,\cdot\}_s,H_s^+,H_s^-)$
is a (not necessarily continuous) family of symplectic splittings.

We recall that the natural space $\bbb_s:= \bigl(D_{\mmax,s}/D_m,\{\cdot,\cdot\}_s\bigr)$
of all boundary values
is a strong symplectic Hilbert space. We assume that there is a family of symplectic splitting
$(\bbb_s,\{\cdot,\cdot\}_s,\bbb_s^+,\bbb_s^-)$ such that
$
H_s^\pm=\bbb_s^\pm\cap D_W/D_m\/.
$

\smallskip

\noi (c) The third group of assumptions concerns the domains $D_s$ and
the reduced Cauchy data:
\begin{enumerate}

\item We are given a continuous family $\{\la_s\}$ of Lagrangian subspaces of $L$.

\item We set
\[
D_s:=\{x\in D_W|(x+D_m)/D_m\in\la_s\}\/.
\]
We assume that $\bigl\{A_{s,D_s}:=A_s^*|_{D_s}; s\in [0,1]\bigr\}$ is a
continuous family (in the gap topology) of self-adjoint Fredholm operators
in $X$.

\item For real $a$ with $\abs{a}\ll 1$, let
\[
\wt\mu_{s,a}:= \bigl(D_m+\ker(A_s^*+aI)|_{D_W}\bigr)/D_m\/.
\]
We assume that there exists a continuous two-parameter family $\{\mu_{s,a}\}$
of Lagrangians in $L$
for $\abs{a}\ll 1$ such that $\mu_{s,a}\,\cap\, \bigl(D_W/D_m\bigr) = \wt\mu_{s,a}$\/.

\item We assume that $(\la_s,\mu_{s,a})\in\Ff\Ll^2(L)$ of index 0.

\item Finally, we make the regularity assumption
\[
\la_s\cap\mu_{s,a}\< D_W/D_m\/.
\]
\end{enumerate}
\end{asss}

\medskip

\begin{rem}\label{r:op-continuity} As a simple example let us
consider the case that the domain of $A^*_s$ is fixed. Let us call it
$D_{\mmax}$\/, and set $D_W:=D_{\mmax}$. By standard norm estimates we
obtain at once
\begin{enumerate}
\item The graph norms are uniformly equivalent with the constants
approaching 1 for $s\to s_0$ for any $s_0\in [0,1]$.

\item The family $\{A_s^*\}_{s\in [0,1]}$ is continuous in $\Cc(X)$,
i.e. with respect to the gap norm.

\item If $\{D_s/D_m\}$ is a continuous family of Lagrangian subspaces of $D_W/D_m$, then
$\bigl\{A_s^*|_{D_s}\bigr\}$ is a continuous family of self-adjoint operators
in $\Cc(X)$.

\end{enumerate}
\end{rem}

\begin{lemma}\label{l:bofu98}
Under Assumptions \ref{a:continuity-symplectic-splitting},
we have for each fixed $s\in [0,1]$ and $\e>0$ sufficiently small
\[
\SF\Bigl\{A_{s,D_s}+aI; {a\in [0,\e]}\Bigr\}
= -\Mas \Bigl\{\la_s,\mu_{s,a}\/;P_s;{a\in [0,\e]}\Bigr\},
\]
where $P_s:L=L_s^+\oplus L_s^-\to L_s^+$ denotes the symplectic splitting projections.
\end{lemma}

\begin{proof} By (a) and c(ii) of Assumptions \ref{a:continuity-symplectic-splitting} and
Corollary \ref{c:stable-ucp}, the operators $A_s+aI$, $a\in [0,\e]$ satisfy weak inner UCP for $\e\ll 1$\/.
From c(ii) and c(iii) of Assumptions \ref{a:continuity-symplectic-splitting}
we have $\la_s\cap\ol{\wt\mu_{s,a}}^{\bbb_s}
= \la_s\cap\wt\mu_{s,a}$\/.
From Assumption \ref{a:continuity-symplectic-splitting}c(ii) it follows that
$D_s/D_m=\la_s\,\cap\, \bigl(D_W/D_m\bigr)$ is a Lagrangian subspace of $\bbb_s$\/,
and so, it makes a Fredholm pair
with the natural Cauchy data space $\g_s(\ker(A^*_s+aI))$ in $\bbb_s$.
Here $\g_s: D_{\mmax,s}\to \bbb_s$ denotes the projection.

Note that $\la_s\,\cap\, \bigl(D_W/D_m\bigr)$ is also
Lagrangian in the weak symplectic Hilbert space $D_W/D_m$ by Lemma \ref{l:embedding}.
A consequence of assumption c(ii) in combination with Proposition
\ref{p:cauchy} is that $(\la_s\cap D_W/D_m,\wt\mu_{s,a})\in\Ff\Ll^2(D_W/D_m)$ and
of index 0. By Theorem \ref{t:cd_continuity}, we have two continuous families
$$\Bigl\{\g_s\bigl(\ker
A_s^*+aI\bigr); {a\in [0,\e]}\Bigr\}\tand \Bigl\{\g_s\bigl(\ker
A_s^*|_{D_W}+aI\bigr); {a\in [0,\e]}\Bigr\}\/.$$

We denote by $\wt P_s:D_W/D_m=H_s^+\oplus H_s^-\to H_s^+$ the symplectic splitting projections
in $D_W/D_m$\/.

From \cite[Theorem 5.1]{BoFu98}, Corollary \ref{c:mas-bofu=mas-zhu}
and Lemma \ref{l:maslov-embedding}, we have
\begin{align*}
\SF&\Bigl\{A_{s,D_s}+aI; {a\in [0,\e]}\Bigr\} \\
&= \Mas_{\operatorname{BF}}
\Bigl\{\g_s(\ker(A^*_s+aI)),
\la_s \cap \bigl(D_W/D_m\bigr);{a\in [0,\e]}\Bigr\} \text{ in $\bbb_s$}\\
&= -\Mas \Bigl\{\la_s \cap \bigl(D_W/D_m\bigr),\g_s(\ker(A^*_s|_{D_W}+aI));
\wt P_s\Bigr\}  \text{ in $D_W/D_m$}\\
&= -\Mas \Bigl\{\la_s,\mu_{s,a}\/;P_s; {a\in [0,\e]}\Bigr\}  \text{ in $L$}.
\end{align*}
\end{proof}

\begin{rem}\label{r:fixm}
In principle, the third line of the preceding formula alignment is the only
place in the proof of our new general spectral flow formula where we really
need our new weak symplectic analysis of Section \ref{s:weak}. As a matter of fact, we could
have avoided weak symplectic analysis by exploiting
the {\it Criss-Cross Theorem} of \cite{BoFuOt01}.
Anyway, the goal is the same, namely to come
\begin{itemize}
\item from more easily
obtainable results in the natural boundary value space $\bbb$
(which, however, is a distribution space in applications)

\item to more easily applicable
results in the familiar $L^2$ space (our $L$ in applications).
\end{itemize}

However, that theorem is quite delicate because it relates curves of Fredholm
pairs of Lagrangians and the Maslov index in two symplectic spaces where
none of them is contained in the other one.

Going via the weak symplectic Hilbert space $D_W/D_m$ has the advantage
that every step can be made by plain embedding arguments, first embedding
$D_W/D_m$ into $\bbb$; then embedding $D_W/D_m$ into $L$.

Actually, the preceding lemma and proof can be considered as an elementary
and intuitive proof of parts of the Criss-Cross Theorem, though not of
the continuity implications which here are assumed and not proved.
\end{rem}

\medskip

Now our main result follows at once.

\begin{theorem} {\rm (General Spectral Flow Formula).}\label{t:gensff}
Under Assumptions \ref{a:continuity-symplectic-splitting} we have
\begin{equation}\label{e:gensff}
\SF\{A_{s,D_s}\}
= -\Mas \{\la_s,\mu_{s,0}\}.
\end{equation}
\end{theorem}

\begin{proof}
To begin with, we sharpen the observation made at the beginning of the proof of the preceding lemma: by (a) and c(ii) of Assumptions \ref{a:continuity-symplectic-splitting} and
Corollary \ref{c:stable-ucp}, for each $s_0$ there exists an $\e(s_0)>0$ such that
the operators $A_s+aI$ satisfy weak inner UCP for all $s,a$ with $\abs{s-s_0},
\abs{a}<\e(s_0)$. Here we use the continuity of the family $\bigl\{A_s^*|_{D_W}\bigr\}$
as bounded operators from $D_W$ to $X$. Since $[0,1]$ is compact, there exists
an $\e>0$ such that the operators $A_s+aI$ satisfy weak inner UCP for all $s\in [0,1]$ and
$\abs{a}<\e$.

We only need to prove the formula \eqref{e:gensff} in a small interval $[s_0,s_1]$.
We consider the two-parameter families
\[
\{A_{s,D_s}+aI\} \tand \{\la_s,\mu_{s,a}\}
\]
for $s\in[s_0,s_1]$ and $a\in[0,\e]$. Because of the homotopy invariance of spectral flow and Maslov
index, both integers must vanish for the boundary loop going counter clockwise around the rectangular domain
from the corner point $(s_0,0)$ via the corner points $(s_1,0)$, $(s_1,\e)$, and $(s_0,\e)$ back
to $(s_0,0)$.

Moreover, for $s_1$ sufficiently close to $s_0$ we can choose $\e$ sufficiently small so
that $\ker(A_{s,D_s}+\e I)=\{0\}$ for all $s\in[s_0,s_1]$. Hence, spectral flow and Maslov index
must vanish on the top segment of our box.

Finally, by the preceding lemma, the left and the right side segments of our curves yield
vanishing sum of spectral flow and Maslov index. So, by additivity under catenation, our
assertion follows.
\end{proof}

By Remark \ref{r:fixm}, we have for fixed maximal domain:

\begin{corollary}\label{c:fixm}
Let $\{A_s:D_m\to X\}_{s\in [0,1]}$ be a family of closed symmetric densely
defined operators satisfying weak inner UCP, with three Hilbert spaces with inclusions
\[
D_m\into D_W\into X,
\]
where $D_m$ is a closed subspace of $D_W$, $D_W=\dom A_s^*$. We assume
that on $D_W$ the graph norms induced by $A_s^*$ and the original norm are
equivalent. Assume that $\{A_s^*:D_W\to X\}$ is a continuous family of
bounded operators. If $\{D_s/D_m\}$ is a continuous family of Lagrangian
subspaces of $D_W/D_m$, such that all
$\bigl(D_s/D_m,\gamma_s(\ker(A_s^*))\bigr)$ are Fredholm pairs in
$D_W/D_m$\/, then $\bigl\{A_{s,D_s}\bigr\}$ is a continuous family of
self-adjoint operators in $\Cc(X)$, and we have
\[
\SF\{A_{s,D_s}\}
= -\Mas \{\g(D_s),\g(\ker A_s^*)\}.
\]
\end{corollary}



\addtocontents{toc}{\medskip\noi}

\section{Curves of Well-posed Boundary Problems}\label{s:geometric}

After having expanded weak symplectic linear algebra and analysis to some length
and detail in the two preceding sections, we shall turn to the geometric
setting and the geometric applications.

\subsection{Symmetric Generalized Dirac Type Operators}\label{ss:gen-di-op}
Let $E$ be a vector bundle over a smooth compact manifold $M$ with
boundary $\Si$\/.
Let $h_1$ be a Hermitian metric on $E$ and $g_1$ a Riemannian metric on $M$. Let
\[
A:\Ci_0(M;E)\to\Ci_0(M;E),
\]
be a linear elliptic differential operator of first order. We assume that
$A$ is formally self-adjoint with respect to the structures $(g,h)$, i.e., it is
symmetric as unbounded densely defined operator in the Hilbert space
$L^2(M;E;g,h)$. For better reading, we shall use the same letter $A$ for the
original operator defined on the space $\Ci_0(M;E)$ of all smooth sections
with support in $M\setminus \Si$, and for its standard extensions to
$\Ci(M;E)$ and to the corresponding first Sobolev spaces $H^1_0(M;E)$
and $H^1(M;E)$.

To give a better understanding of Assumption (ii) in Definition
\ref{d:dirac} below, we shall describe what is meant by product
form of an operator near the boundary - without assuming product
metric structures from the very beginning. The details are worked
out separately, see M. Lesch \cite{Le04a}.

We parametrize a collar neighbourhood $N\< M$ of $\Si$ in $M$ by
a diffeomorphism (a {\em partial isometry})
\[\begin{matrix}
\f&:&N&\to &[0,\e)\times \Si\\
\ &\,&p&\mapsto &(t,\eta)
\end{matrix}
\]
with $\f|_{\Si}=\id|_{\Si}$\,.
Here $t$ denotes the inward normal coordinate in $N$. We choose
a metric $g$ on $M$ such that
\begin{equation}\label{e:induced-metric}
\f_*g|_N= dt^2\wedge g_{\Si} \text{ with $g_{\Si}:=g_1|_{\Si}$}\,.
\end{equation}
We denote the retraction of $N$ onto $\Si$ by $\rho$. Let
\[
\F:L^2\bigl(N;E|_N\bigr)\too L^2\bigl([0,\e)\times \Si;\rho^*(E|_{\Si})  \bigr)
\]
be a unitary operator covering $\f$. Then $A$ can be written on $N$ in
the form
\begin{equation}\label{e:generalized_dirac_type}
\F A|_N \F\ii = J_t(\frac{\dpa}{\dpa t}+B_t) \,,
\end{equation}
where $\{J_t\in \Hom(E|_{\Si})\}$ is a curve of bundle isomorphisms and $\{B_t\}$ a curve
of linear elliptic differential operators in `tangential'
direction; i.e., on $\Si$.
Since $A$ is formally self-adjoint, it follows
that $J_t^*=-J_t$ and $\frac{\dpa}{\dpa t}\bigl(J_t\bigr)-B_t^*J_t=J_tB_t$\,.

\medskip

\begin{definition}\label{d:dirac}
We shall call such an operator $A$ a {\em symmetric generalized Dirac type operator},
if the following conditions are satisfied.

\begin{enumerate}

\item We require that each connected component $M$ has non-empty boundary.

\item The principal symbol of $B_0$ is symmetric, i.e., $B_0^*-B_0$ operator of order 0.

\item Moreover, we shall assume that $A$ satisfies the weak inner unique continuation
property (UCP) with respect to $\Si$, i.e.,
$\ker A^*\,\cap\, H_0^1(M;E)= \{0\}$, where $A^*$ denotes the maximal closed extension of $A$
and $H_0^1(M;E)$ denotes the subspace of the first Sobolev space $H^1(M;E)$ comprising elements with
support in the interior of $M$.
\end{enumerate}
\end{definition}

\begin{rem}\label{r:generalized-dirac}
(a) It is well known that each compatible Dirac operator (in the true geometric sense,
i.e., given by composition of Clifford multiplication with a compatible connection
in the underlying bundle
of Clifford modules) is a generalized Dirac type operator in the
sense of Definition \ref{d:dirac} (see also \cite[Lemma 2.2]{BoMaWa02} or G. Grubb
and R. T. Seeley\cite{GrSe95}).

\noi (b) The point of assumption (ii) of the preceding definition is that we only assume
symmetric principal symbol for the tangential operator and only at the precise boundary.
Compared with the usual assumptions on Dirac type operators there are two
additional {\it generalizations}:
that we admit variable $J_t$ and that we do no longer
assume that the symbol $J_0$ is unitary. In particular, we drop
the common assumption $J_0^2=-I$ which comes from Clifford multiplication.
Motivated mainly by applications in the field of Hamiltonian dynamics, we assume only
the invertibility of $J_0$ and $J_t^*=-J_t$ for small $t$.
\end{rem}

\medskip

One way to address parameter dependence is to consider a (big) Hermitian vector bundle $\EE$ over a (big)
compact Hausdorff space $\MM$. We assume that $\MM$ itself is a fibre bundle over the interval
$[0,1]$ such that $\MM$ is a continuous family of compact smooth Riemannian
manifolds $j_s:M_s \into \MM$ with boundary $\Si_s$.
We require that the vector bundle structure is compatible with the boundary part.
More precisely, we shall have a trivialization
\[
\f:M_0\times [0,1]\simeq \MM
\]
such that $\pi\circ\f\ii\circ j_s: (M_s,\Si_s)\to (M_0,\Si_0)$ is a diffeomorphism. Here $\pi$
denotes the natural projection $\pi:M_0\times [0,1]\to M_0$\,. We do not assume that
$M_s$ or $\Si_s$ are connected. Note that the trivializations define smooth structures on
$\ran\f|_{M_0\setminus \Si_0\,\times\, [0,1]}$ and so on $\MM$ and $\EE$.

We consider a smooth linear differential
operator $\AA:\Ci(\MM;\EE)\to\Ci(\MM;\EE)$ which induces a smooth family of generalized Dirac type operators
over $M_s$
\begin{equation}\label{e:as-original}
A_s:\Ci_0(M_s;E_s) \to \Ci(M_s;E_s),
\end{equation}
where $E_s\to M_s$ denotes the induced bundle, i.e., the pull back $j_s^*(\EE)$.
Recall that $\Ci_0(M_s;E_s)$ denotes the space of smooth sections with
support in the interior $M_s^0:=M_s\setminus \Si_s$ of $M_s$\,. We define the Hilbert space
\[
  H^1_0(M_s;E_s):= \ol{\Ci_0(M_s;E_s)}^{H^1(M_s;E_s)}\,,
\]
where $H^1(M_s;E_s)$ denotes the first Sobolev space. The inner product is given by
the graph inner product. The corresponding graph norm is, again by ellipticity (see
the relevant estimates, e.g., in \cite[Chapter 18]{BoWo93}) equivalent to the Sobolev norm.
Then the operator $A_s$ extends to a bounded operator
\begin{equation}\label{e:as-h1}
A_s:H^1(M_s;E_s) \to L^2(M_s;E_s),
\end{equation}
which we denote by the same symbol.

Notice that we have a continuous family of bundle isomorphisms $T_s:E_s\to E_{s_0}$ for
$|s-s_0|\ll 1$ induced by a local trivialization of the fibre bundle $M\to[0,1]$.
By construction the induced family
\begin{equation}\label{e:operator_comparison}
\bigl\{T_s\circ A_s\circ T_s\ii
: H^1(M_{s_0};E_{s_0}) \to L^2(M_{s_0};E_{s_0})\bigr\}
\end{equation}
of bounded operators is continuous.

Since the operator $A_s$ is elliptic, the extension to $H^1_0(M_s;E_s)$
is a closed operator, if we view it as an unbounded operator in
$L^2(M_s;E_s)$ (see, e.g., \cite[Proposition 20.7]{BoWo93}).

We assume that
all $A_s$ are formally self-adjoint, i.e., $A_s\< A_s^*$\,.
Note that we do not make any assumptions about product structures near the boundary $\Si_s$\,.

\medskip

For each $s$ we choose a well-posed self-adjoint boundary condition
$P_s\in\Grass_{\sa}(A_s)$ in the sense of Br{\"u}ning and Lesch \cite{BrLe01}.
That is a pseudo-differential projection
$P_s: L^2(\Si_s;E_s|_{\Si_s}) \to L^2(\Si_s;E_s|_{\Si_s})$ which defines a
self-adjoint Fredholm extension $A_{s,P_s}$ in $L^2(M_s;E_s)$ by
\[
\dom A_{s,P_s} = D_s := \{x\in H^1(M_s;E_s)\mid P_s(x|_{\Si_s})=0\}.
\]

To define the continuity of a family $\{P_s\}_{s\in [0,1]}$ of boundary conditions at $s_0\in [0,1]$,
we use again the family $T_s:E_s\to E_{s_0}$\,. We shall call the family
$\{P_s\}$ continuous if and only if the family
\[
\bigl\{T_s|_{({E_s}|_{\Si_s})}\circ P_s\circ T_s\ii|_{({E_{s_0}}|_{\Si_{s_0}})}
: L^2(\Si_{s_0};E_{s_0}|_{\Si_{s_0}}) \to L^2(\Si_{s_0};E_{s_0}|_{\Si_{s_0}})\bigr\}
\]
of bounded operators in $L^2(\Si_{s_0};E_{s_0}|_{\Si_{s_0}})$ is continuous.
It follows that $\{T_s\circ A_{s,P_s}\circ T_s\ii\}$ is a continuous curve in the space
$\Cc\Ff(L^2(M_{s_0};E_{s_0}))$ of closed Fredholm operators in the fixed
Hilbert space $L^2(M_{s_0};E_{s_0})$. In general, the transformed operators
are no longer self-adjoint. However, the transformation does not change the
spectrum and it turns out that the integer $\SF\{T_s A_{s,P_s} T_s\ii\}$ is
well-defined and does not depend on the choice of
$T_s$\,.

Geometrically speaking, we do not deal with continuous curves in a
fixed space of self-adjoint Fredholm operators but with sections of a
homomorphism bundle of self-adjoint Fredholm operators between Hilbert space bundles
over the interval $[0,1]$. In that way, the spectral flow $\SF\{A_{s,P_s}\}$ is defined.

\medskip

Let $Q_s: L^2(\Si_s;E_s|_{\Si_s}) \to L^2(\Si_s;E_s|_{\Si_s})$ denote the
Calder{\'o}n projection belonging to the operator $A_s$. It is a pseudo-differential projection with
\[
\ran Q_s = \ol{\{x|_{\Si_s} \mid x\in \Ci(M_s;E_s) \tand A_s^*x=0
\text{ in $M_s\setminus \Si_s$}\}}^{L^2(\Si_s;E_s|_{\Si_s})}
\]
the $L^2$--Cauchy data space of $A_s$\,. From the assumed weak inner
Unique Continuation Property (UCP) $\ker A_s=\{0\}$
(here $A_s$ denotes the closed extension of \eqref{e:as-h1} and not the original
operator of \eqref{e:as-original})
we can deduce that the Calder{\'o}n
projections are continuous in this sense, i.e., that for $|s-s_0|\ll 1$
the family $\{T_s|_{({E_s}|_{\Si_s})}\circ Q_s\circ T_s\ii|_{({E_s}|_{\Si_s})}\}$
is continuous. This was shown in \cite[Section 3.3]{BoFu98} in a broad functional analytical
setting, but under the
narrow assumption that $A_s$ is just a perturbation of a fixed closed symmetric
operator $A_0$ by a continuous family of bounded self-adjoint operators, i.e., in
our applications only variation of the zero'th order coefficients are admitted;
\cite[Theorem 3.9]{BoLePh01} permits also variation of the first order coefficients,
but only for {\it classical} Dirac type operators, i.e., for unitary $J_{t,s}$ with
$J_{t,s}^2=-I$; a proof for {\it generalized} Dirac type operators in the sense of
Definition \ref{d:dirac} follows from \cite{Le04a}.

Then the family $\{(\ker P_s,\ran Q_s)\}$ (we suppress the local transformations $T_s$) is a
continuous curve of Fredholm pairs of closed subspaces of $L^2(\Si_s;E|_{\Si_s})$ (all of index 0).
Now we have two invariants. On one side, as mentioned before, we have the spectral flow
$\SF\{A_{s,P_s}\}$. It is a {\em spectral} invariant, defined by the behaviour of the spectrum
near the 0-eigenvalue as the parameter $s$ runs from 0 to 1.
On the other side, we have these pairs of
$L^2$-sections over $\Si_s$\,. Under suitable assumptions, they define a second invariant,
the Maslov index of these pairs. It is a {\em classical} invariant of symplectic analysis.
Roughly speaking, it is given by
counting how many times (and with which sign and which multiplicity) the one family, $\{\ker P_s\}$
comes to share a non-trivial subspace with the other family $\{\ran Q_s\}$
when $s$ runs from 0 to 1.
In applications, the Maslov index is easier to determine than the spectral flow.

We shall express the spectral flow of the operator family $\{A_{s,P_s}\}_{s\in [0,1]}$
as the Maslov index of the curve of pairs $\{(\ker P_s, \ran Q_s)\}_{s\in [0,1]}$.
To do that, however, we must overcome some rather serious difficulties. The main reason for our difficulties
is that, in general, the {\it natural} boundary value spaces $\beta_s:= D_{\mmax}(A_s)/D_m(A_s)$
vary in an uncontrollable way, both as Hilbert spaces and as symplectic spaces. Recall that
$D_m(A_s)= H_0^1(M_s;E_s)$ denotes the (minimal) domain of the closed symmetric operator $A_s$
which is fixed or easily fixable. On the contrary, the space $D_{\mmax}(A_s)$ which
denotes the domain of the maximal closed extension of $A_s$ may vary in an uncontrollable way.
Roughly speaking, that is the reason why we cannot easily follow the functional
analytical path worked out in \cite{BoFu98}.

\medskip

Instead of working in the natural distribution space $\beta_s$ we must consider
the space of {\it reduced} boundary values
\begin{equation}\label{e:h_one_half}
H^{1/2}(\Si_s;{E_s}|_{\Si_s})= \g\bigl(H^1(M_s;E_s)\bigr)
\simeq H^1(M_s;E_s)/H^1_0(M_s;E_s)
\end{equation}
with trace map $\g: H^1(M_s;E_s) \to
H^{1/2}(\Si_s;{E_s}|_{\Si_s})$, coming from the geometrically fixable first Sobolev space.
For each $s$ we have a naturally defined
symplectic structure on $H^{1/2}(\Si_s;{E_s}|_{\Si_s})$ given by Green's form
\[
\w_s(\g(x),\g(y)) :=  \lla A^*_s x,y\rra_{L^2(M_s;E_s)} - \lla x, A^*_s y\rra_{L^2(M_s;E_s)}
\]
for $x,y\in H^1(M_s;E_s)$.

\smallskip

The first difficulty is that these symplectic
structures are varying when $s$ varies. Note that $\w_s$ extends to
\[
\w_s: L^2(\Si_s;E|_{\Si_s}) \times L^2(\Si_s;E|_{\Si_s}) \to \C\,.
\]
Applying the bundle isomorphisms $\{T_s\}$ close to one $s_0$, we fix the Hilbert space
$L^2(\Si_{s_0};E|_{\Si_{s_0}})$ -- with varying symplectic structures $(T_s^{-1})^*\w_s$.
This is the main technical problem we shall handle in this paper.

\medskip

For {\it fixed} symplectic structures one has a nice reduction method to obtain a general
spectral flow formula for
the spectral flow $\SF\{A_{s,P_s}\}$ of a curve of varying operators with varying domains.
Then, roughly speaking, the problem
can be reduced to calculating the spectral flow for two curves $\{A_{0,P_{0,t}}\}$ and
$\{A_{1,P_{1,t}}\}$. Each of the two curves has a fixed operator and solely varying domain.
This reduction can be achieved by uniformly connecting
each operator $A_{s,P_{s}}$ to an invertible operator $A_{s,Q_{s}}$ within the category
of well-posed self-adjoint boundary value problems (see, e.g., \cite{KiLe00} where, however, advanced
results about the $\eta$-invariant are used to get further).

This approach is not open in our situation where we admit varying symplectic structures and,
hence, are prevented from the deformation and reduction approach.

\smallskip

Having fixed the space $H^{1/2}(\Si_s;E|_{\Si_s})$ as
the natural symplectic Hilbert space for traces at the boundary in our setting, another delicate problem
arises: our proof of Theorem \ref{t:gsff1} will show that, in general, one can not expect
the validity of the spectral flow formula $\SF\{\dots\} = -\Mas\{(\dots,\dots)\}$ when taking
the Maslov index on the right in the (strong) symplectic space $L^2$\,.
Roughly speaking, the intersection dimensions entering into the Maslov index
on the right side of the formula can become too large in $L^2$\,. It is a special feature
of pseudo-differential projections that any Fredholm pair $(\ker P_s,\ran Q_s)$ of
Lagrangian subspaces of $L^2(\Si_s;E|_{\Si_s})$ has its intersection
\begin{equation}\label{e:regularity}
\ker P_s\, \cap\, \ran Q_s\, \< \, H^{1/2}(\Si_s;E|_{\Si_s})
\end{equation}
for all $s$ by a kind of elliptic regularity.

Here it is important that the symplectic Hilbert space structures of $H^{1/2}(\Si_s;E|_{\Si_s})$
and $L^2(\Si_s;E|_{\Si_s})$ are compatible and their symplectic splittings are defined by
the bundle endomorphisms $J_{s,0}:E|_{\Si_s}\to E|_{\Si_s}$ in the following way:
\begin{multline}\label{e:hpm_in_h12}
H_s^\pm := H^{1/2}(\Si;E_s^\pm|_{\Si_s}) \tand L_s^\pm :=
L^2(\Si;E_s^\pm|_{\Si_s})
\\ \text{ with }
E_s^\pm|_{\Si_s}:= \text{ lin. span
of}\left\{\begin{array}{l}\text{positive}\\ \text{negative} \end{array}
\right\} \text{ eigenspaces of $iJ_{s,0}$}.
\end{multline}
Note that $L_s^+, L_s^-$ change continuously if $J_{s,0}$ changes continuously.

So, as a
result of our analysis it turns out that the spectral
flow formula {\em is} valid also in $L^2$ for well-posed boundary value problems which are defined
by pseudo-differential projections.

The difference between the spaces $H^{1/2}(\Si_s;E|_{\Si_s})$ and $L^2(\Si_s;E|_{\Si_s})$ must be
emphasized also under the perspective of continuity. In difference to \cite{BoFu98} where
$D_{\mmax}(A_s)$ and the Fredholm domain (say, given by $P_s=P_0$) are kept fixed, now
the continuity of the operator family $\{A_{s,P_s}\}$ is a problem.
It can be more easily derived from the
continuity of $\{P_s\}$ (or $\{\ker P_s\}$) in $H^{1/2}$ than in $L^2$\,.
However, our pair $(\ker P_s,\ran Q_s)$ is, a priori, only continuous in $L^2(\Si_s;E|_{\Si_s})$.
In general, $\ker P_s\,\cap\, H^{1/2}(\Si_s;E|_{\Si_s})$ is not continuous
in $H^{1/2}(\Si_s;E|_{\Si_s})$ without making additional assumptions.

\medskip

The price of working in $H^{1/2}(\Si_s;E|_{\Si_s})$, or, more precisely
(after locally transforming)
in $H^{1/2}(\Si_{s_0};E|_{\Si_{s_0}})$ must be paid by loosing the invertibility of the
associated structures $J'_s$ which are defined by
\[
(T_s\ii)^*\w_s(x,y)
= -\lla J'_s x,y\rra_{H^{1/2}(\Si_{s_0};E_{s_0}|_{\Si_{s_0}})} \text{ for $x,y \in
H^{1/2}(\Si_{s_0};E_{s_0}|_{\Si_{s_0}})$}.
\]
So, we have to deal with \textit{weak} symplectic structures as treated in Section \ref{s:weak}.

\smallskip

\subsection{Proof of Theorem \ref{t:gsff1} and Theorem
  \ref{t:split} }

We fix the manifolds and bundles and drop the fibre bundle notation and the
transformations for easier reading.

\begin{rem}\label{r:weak-to-inner}
To prove Theorem \ref{t:gsff1} and Theorem \ref{t:split}, we can weaken our assumption
regarding weak inner UCP. If the principal symbol of the tangential operators $B_t$
is symmetric in a collar of $\Si$, we can show that weak UCP with respect to $\Si$
in the smooth category (i.e., $\ker A\cap \Ci_0(M;E)=\{0\}$) implies weak inner UCP in the sense of
Definition \ref{d:ucp_general} (i.e., $\ker A^*\cap H_0^1(M;E)=\{0\}$).

The argument runs as follows:
Let $x\in \ker\bigl(A:H_0^1(M;E)\to L^2(M;E)\bigr)$. Since $A$ is
elliptic, $x$ is smooth up to the boundary. Assuming the symmetry of the principal
symbol of the tangential operators in the collar we can find a continuation of $A$
to an operator $A'$ in a manifold $M'=[-\d,0]\cup_{\Si}M$ of the same type.
Repeating the arguments of the
proof of \cite[Theorem 2.7]{BoMaWa02} we obtain that the support of $x$ is
totally contained in $M\setminus \Si$. The details will be worked out separately.
So, by our weakened assumption, we have $x=0$.
\end{rem}

\medskip

\begin{proof}[Proof of Theorem \ref{t:gsff1}]
Set
\begin{align*}
&L:=L^2(\Si;E|_{\Si}),\, D_W:=H^1(M;E),\\
& D_m:=H_0^1(M;E),\,
\bbb_s:=(\dom A^*_s)/D_m
\end{align*}
with the symplectic forms induced by
\begin{equation}\label{e:dirac-omega}
\{x,y\}_s := -\int_{\Si} \lla J_{s,0}x,y\rra d\Si\/.
\end{equation}
The inclusions of Assumption \ref{a:continuity-symplectic-splitting}a
are clear. The splittings
$L_s^\pm$ and $H_s^\pm$ are defined in \eqref{e:hpm_in_h12};
we set $\bbb_s^\pm:=\ol{H_s^\pm}^{\bbb_s}$\/.

By Br{\"u}ning and Lesch \cite[Theorem 1.5]{BrLe01}, the operators
$A_{s,P_s}$ are self-adjoint Fredholm operators. By \cite[Theorem
3.9b]{BoLePh01} (an elaboration of the arguments will be published
separately), the projections $Q_s$ are pseudo-differential
projections and $\bigl\{\range Q_s\bigr\}$ is continuous in
$L^2(\Si;E|_{\Si})$. So, $\Mas\{\ker P_s,\range Q_s\}$ is well
defined.

By \cite[Theorem 3.9d]{BoLePh01}, family $\{A_{s,P_s}\}$ is continuous in
$\Cc(L^2(M;E))$. So,  $\SF\{A_{s,P_s}\}$ is well defined.

Regularity \eqref{e:regularity} is proved like in \cite[Chapter 19]{BoWo93}
or \cite[Theorem 6.5]{BrLe01}. By our Theorem \ref{t:gensff}, we have the assertion.
\end{proof}

\medskip

\begin{proof}[Proof of Theorem \ref{t:split}]
Let $M^{\sharp}$ denote the compact manifold
$M^+\sqcup M^-=
\bigl(M\setminus \Si\bigr)\,\cup\, \bigl((\Si\sqcup(-\Si)\bigr)$ with boundary
$\dpa M^{\sharp} =\dpa M^+\sqcup\dpa M^-=\Si \sqcup(-\Si)=:\Si^{\sharp}$
and $E^{\sharp}\to M^{\sharp}$ the corresponding Hermitian bundle.
Fixing $\Si$ induces a
decomposition
\[
L^2(M;E)\cong L^2(M^+;E|_{M^+})\oplus L^2(M^-;E|_{M^-})
=L^2(M^{\sharp};E^{\sharp}),
\]
and for the first Sobolev space
\[
H^1(M^+;E|_{M^+})\oplus H^1(M^-;E|_{M^-})
=H^1(M^{\sharp};E^{\sharp}).
\]
Correspondingly we obtain an operator $A_s^{\sharp}$ for each $s\in [0,1]$
which is a symmetric generalized operator of Dirac type according to the
assumptions made for Theorem \ref{t:split}.

For the Calder{\'o}n projection of $A^{\sharp}$ we have
\begin{align*}
\range Q_s^{\sharp}&=\range Q_s^+\boxplus \range Q_s^-\\
&\< L^2(\Si^{\sharp}; E^{\sharp}|_{\Si^{\sharp}})
= L^2(\Si; E|_{\Si})\boxplus L^2(-\Si; E|_{-\Si}).
\end{align*}
Here $L^2(\Si^{\sharp}; E^{\sharp}|_{\Si^{\sharp}})$ is provided with the symplectic
forms $\w_s^{\sharp}:=\w_s\boxplus(-\w_s)$ where $\w_s$ is defined like in
\eqref{e:dirac-omega}.

Let $\D$ denote the diagonal in $L^2(\Si; E|_{\Si})
\boxplus L^2(-\Si; E|_{-\Si})$. Clearly, for all $s$
it is a Lagrangian subspace with respect to $\w_s^{\sharp}$ and makes a Fredholm pair
with each $\range Q_s^{\sharp}$\/. The projection of
$L^2(\Si^{\sharp}; E^{\sharp}|_{\Si^{\sharp}})$ onto $\D$
satisfies the regularity condition \eqref{e:regularity}
(even it is not a pseudo-differential operator
over the manifold $\Si^{\sharp}$, as noticed in \cite[Section 5]{KiLe00}).

Consequently, we have on the manifold $M^{\sharp}$
a natural self--adjoint elliptic boundary condition
(in the sense of our Theorem \ref{t:gsff1}) defined for $A_s^{\sharp}$
by the {\it pasting} domain
\begin{align}\label{e:pasting}
D^\sharp :&=\{(x,y)\in H^1(M^{\sharp};E^{\sharp})  \mid \g^+(x)=\g^-(y)\}\\
&= \{(x,y)\in H^1(M^{\sharp};E^{\sharp})  \mid \g^{\sharp}(x,y)\in\D\},
\end{align}
where $\g^\pm:H^1(M^\pm;E|_{M^\pm})\to H^{\12}(\pm\Si;E|_{\pm\Si})$
and $\g^{\sharp}:H^1(M^{\sharp};E^{\sharp})\to H^{\12}(\Si^{\sharp};
E^{\sharp}|_{\Si^{\sharp}})$ denotes the trace maps.
Let $A^{\sharp}_{s,\D}$ denote the operator which acts like $A^{\sharp}_s$ and has domain $D^{\sharp}$\/.

By these definitions and applying Proposition \ref{p:boxplus}b
(switch of symplectic forms) and Theorem \ref{t:gsff1} to the operator family
$\{A^{\sharp}_{s,\D}\}$ we obtain
\begin{align*}
\SF\{A_s\}&= \SF\{A^{\sharp}_{s,\D}\}
\fequal{Th. 0.1}-\Mas\{\D,\range Q_s^+\boxplus \range Q_s^-\}\\
&\fequal{\eqref{e:maslov4}} -\Mas\{\range Q_s^-, \range Q_s^+\} \
\quad\quad \text{in}\;L^2(\Si; E|_{\Si})\\
&\fequal{\eqref{e:maslov3}}  \Mas\{\range Q_s^+, \range Q_s^-\} \
\quad\qquad \text{in}\;L^2(-\Si; E|_{-\Si})\\
&\fequal{Th. 0.1}  \SF\{A_{s,I-Q_s^+}^-\}.
\end{align*}\end{proof}

\begin{note}
If one is only interested in the equality of the spectral flows on the whole
manifold and on the part, on needs not to argue with the Maslov index, as we do, but
can find a direct proof in \cite[Corollary 5.6]{KiLe00}
  based solely on the homotopy
invariance of the spectral flow of a related two-parameter family.
\end{note}

\smallskip

\subsection{Applications to Hamiltonian Dynamics}\label{ss:hamiltonian}

\begin{proof}[Proof of Theorem \ref{t:First_Order}]
This theorem extends the validity of \cite[Proposition 3.2]{Zh01}. The difference is,
that now we admit that the matrices $j_{s,t}$ vary with $s$. We can do that, because
we admit varying symplectic structures on the spaces of boundary values.

We have the following basic spaces:
\[
D_m:=H_0^1([0,T],\C^m) \ \tand\ D_{\mmax}=D_W:=H^1([0,T],\C^m).
\]
By Sobolev embedding theorem, we have
\[
H^1([0,T];\C^m)\< C([0,T];\C^m).
\]
So, we can describe explicitly
\[
\begin{matrix}
\g&:& H^1([0,T];\C^m) &\too & \bbb_s\\
\ &\ & x &\mapsto & \bigl(x(0),x(T)\bigr)
\end{matrix}
\]
and $\bbb_s:=(\C^m\oplus \C^m,\{\cdot,\cdot\}_s)$ with symplectic
form induced by Green's form
\begin{align*}
\{\g(x),\g(y)\}_s:&=\lla A_s^*x,y\rra - \lla x,A_s^*y\rra\\
&= -\lla j_{s,t}x(t),y(t)\rra\bigr|_0^T
=-\w_s(\g(x),\g(y)).
\end{align*}
By our assumption, $W_s$ is Lagrangian in $\bbb_s$\/, so the operator
$A_{s,W_s}$ (as defined in the introduction) is self-adjoint.
It satisfies weak inner UCP (i.e., $\ker A_s|_{D_m}=\{0\}$)
because of the uniqueness of the solutions for our linear system with
continuous coefficients.
Clearly, $\dim \ker A_{s,W_s}\leq m$. To see that $A_{s,W_s}$ is a
Fredholm operator, we must only check that the range is closed.
This follows from elliptic estimate in the usual way.

By Corollary \ref{c:fixm} we have
\begin{align*}
\SF\{A_{s,W_s}\} &= -\Mas\{\g(\dom A_{s,W_s}),\g(\ker A_s^*)\}\\
&= -\Mas\{W_s,\GGG(\G_s(T))\} \quad \text{ in $\bbb_s$}\\
&= \Mas\{\GGG(\G_s(T)),W_s\} \qquad \text{ in $(\C^m\oplus \C^m,\w_s)$}
\end{align*}

\end{proof}

\smallskip

\begin{proof}[Proof of Theorem \ref{t:hamiltonian}]
We see that
\[
\lla L_sx,y\rra= \lla x,L_s y\rra \ \text{ for all $x,y\in D_m$}\/,
\]
where
\[
D_m:=H_0^2([0,T],\C^m) \ \tand\ D_{\mmax}=D_W:=H^2([0,T],\C^m).
\]
So, $L_s$ is formally self-adjoint. By Sobolev embedding theorem, we have
$H^2([0,T];\C^m)\< C^1([0,T];\C^m)$. Clearly, $D_{\mmax}/D_m=\C^{4m}$ and
\begin{equation}\label{e:trace-second-order}
\g_s(x)=\bigl(u_s(x)(0), u_s(x)(T)\bigr)
=\bigl(u_1(0),u_2(0),u_1(T), u_2(T)\bigr),
\end{equation}
with
\[
\bigl(u_s(x)\bigr)(t):= \bigl(p_s(t) \dd tx(t)+q_s(t)x(t)\,,\,x(t)\bigr)
=:(u_1(t),u_2(t)).
\]
This provides the following symplectic form, induced by Green's form:
\begin{align*}
\{\g_s(x),\g_s(y)\}_s:&= \lla L_sx,y\rra - \lla x,L_s y\rra \\
&= -\lla p_s\dd t x+q_sx,y\rra\bigr|_0^T + \lla x,p_s\dd t y+q_sy\rra\bigr|_0^T\\
&= - \lla u_1(T),v_2(T)\rra +\lla u_1(0), v_2(0)\rra \\
&\qquad + \lla u_2(T),v_1(T)\rra -\lla u_1(0), v_1(0)\rra ,
\end{align*}
where $\g_s(x)$ as in \eqref{e:trace-second-order} and
$\g_s(y)=\bigl(v_1(0),v_2(0),v_1(T), v_2(T)\bigr)$ similarly defined.
Comparing with
\[
\w_s(\g_s(x),\g_s(y))=\Bigl{\lla}
\begin{pmatrix}
-J&0\\
0&J
\end{pmatrix}
\begin{pmatrix}
u_1(0)\\ u_2(0)\\ u_1(T)\\ u_2(T)
\end{pmatrix}\,,\,
\begin{pmatrix}
v_1(0)\\ v_2(0)\\ v_1(T)\\ v_2(T)
\end{pmatrix}
\Bigr{\rra}
\]
yields $\w_s(\g_s(x),\g_s(y))=-\{\g_s(x),\g_s(y)\}_s$\/. So, we have
$\bbb_s= (\C^{4m},-\w_s)$. From the Lagrangian property of $W_s$, we
obtain that the operator $L_{s,W_s}$ is self-adjoint, and from elliptic
estimate that its range is closed and so that it is a Fredholm operator.

To determine the Cauchy data space we notice
\[
x\in\ker L_s\ifff \bigl(u_s(x)\bigr)(t)=\G_s(t)\bigl(u_s(x)\bigr)(0)
\text{ for all $t\in [0,T]$}.
\]
This gives $\g(\ker L_s)=\GGG(\G_s(T))$. Then by Theorem \ref{t:gensff}
and Proposition \ref{p:boxplus} (symplectic inversion) we finally obtain
\[
\begin{matrix}
\SF\{L_s,W_s\} &=& - \Mas\{W_s,\GGG(\G_s(T))\}&\qquad \text{in $(\bbb_s,-\w_s)$}\\
\ &=& \Mas\{\GGG(\G_s(T)),W_s\}&\qquad \text{in $(\C^{4m},\w_s)$}.
\end{matrix}
\]
\end{proof}

\medskip

As a corollary for fixed domain, we consider the index form
\begin{align*}
\Ii_{s,R}(x,y):&=
\int_0^T\bigl(\lla p_{s,t}\dd tx+q_{s,t}x,\dd ty\rra+\lla q_{s,t}^*\dd tx,y\rra
+\lla r_{s,t}x,y\rra\bigr)dt\\
&\text{for $x,y\in H^1([0,T];\C^m)$ with $(x(0),x(T)),(y(0),y(T)) \in R$}.
\end{align*}
Here the coefficients are of the operator $L_s$ of \eqref{e:l-s}, and $R$ denotes
an arbitrary subspace of $\C^{2m}$\/.

We fix our domain by setting
\[
W_s:=W(R)=\{(x,y,z,u)\mid (x,-z)\in R^\perp, (y,u)\in R\}
\]
independently of $s$. Then we have

\begin{corollary}\label{c:morse-index}
\begin{equation}\label{e:gmit}
\SF\{\Ii_{s,R}\}=\SF\{L_{s,W(R)}\}=\iii_{W(R)}(\{\G_1(t)\})
-\iii_{W(R)}(\{\G_0(t)\}),
\end{equation}
where
\[
\iii_{W}(\{\G(t)\}):=\Mas\{\GGG(\G(t)),W\}
\]
denotes the {\em Maslov-Long index} of the symplectic path
$\{\G(t)\}_{t\in [0,1]}$ with respect to a fixed Lagrangian
domain $W \<\C^{4m}$.
\end{corollary}

\begin{note} (a) The Maslov-Long index is more general than the usual
Conley-Zehnder index because $\G(0)\neq I_{2m}$ and general
boundary condition $W\neq\GGG(I_{2m})$ are admitted, see also Long
and Zhu \cite{LoZh00}. We have the relationship
\[
\iii_{\GGG(I_{2m})}(\{\G(t)\})=\iii_{\operatorname{CZ}}(\{\G(t)\})+m
\]
for real path $\{\G(t)\}$ with $\G(0)=I_{2m}$.

\noi (b) If $p_{s,t}$ are positive definite, then the first equality of \ref{e:gmit}
is trivial.
\end{note}

\begin{proof}
By \cite[Theorem 1.1]{Zh01}, we need only to prove the second equality of \ref{e:gmit}.
First, we apply Corollary \ref{c:fixm}. Then we consider the two-parameter
family $\{\G_{s,t}\}_{s\in [0,1],t\in [0,T]}$ and apply homotopy invariance
and catenation additivity of the Maslov index. Since $\G_s(0)=I$, we obtain
\begin{align*}
\SF\{L_{s,{W(R)}}\}&=\Mas\{\GGG(\G_s(T)),{W(R)}\} = \iii_{{W(R)}}(\{\G_s(T)\})\\
&=\iii_{W(R)}(\{\G_s(0)\})+\iii_{W(R)}(\{\G_1(t)\})-\iii_{W(R)}(\{\G_0(t)\})\\
&=\iii_{W(R)}(\{\G_1(t)\})-\iii_{W(R)}(\{\G_0(t)\}).
\end{align*}
\end{proof}


\bigskip

\appendix

\setcounter{secnumdepth}{2}

\setcounter{theorem}{0}

\setcounter{equation}{0}



\addtocontents{toc}{\medskip\noi}
\section{Spectral Flow}\label{s:appendix1}

The spectral flow for a one parameter family of linear
self-adjoint Fredholm operators was introduced by M. Atiyah, V.
Patodi, and I. Singer \cite{AtPaSi75} in their study of index
theory on manifolds with boundary. Since then other significant
applications have been found. Later this notion was made rigorous
for curves of bounded self-adjoint Fredholm operators in J.
Phillips \cite{Ph96} and for continuous curves of self-adjoint
(generally unbounded) Fredholm operators in Hilbert space in
\cite{BoLePh01} by Cayley transform. The notion was generalized to
higher dimensional case in X. Dai and W. Zhang \cite{DaZh98}, and
to more general operators in \cite{Wo85,Zh01,ZhLo99}.

In this Appendix we shall provide a rigorous definition of the {\it
spectral flow} of {\it spectral-continuous} curves of {\it admissible} closed operators
in Banach space relative to a co-oriented real curve  $\ell\<\C$.
(All the preceding terms will be explained).

\begin{note}
Of course, if $\ell$ is an open segment of a straight line, then the following
choices are immediate.
\end{note}

Let $X$ be a Banach space, and $A\in\Cc(X)$. Let $N\<\C$ be a bounded open subset.
Assume that $\s(A)\cap\partial N$ is a finite set.
Then there exists an open subset $\wt N\subset N$ such that
\begin{equation}\label{e:a-rev-admissible}
\ol{\wt N}\< N,\
\partial \wt N \in C^1\/,\  \s(A)\cap \wt N=\s(A)\cap N,
\tand \s(A)\cap\partial \wt N=\emptyset,
\end{equation}
and the spectral projection
\begin{equation} \label{P_N}
P_N(A) := -\frac 1{2\pi i}\int_{\partial \wt N}(A-\la I)^{-1} d\la
\end{equation}
is well defined and does not depend of the choice of $\wt N$.
By Theorem III.6.17 of \cite{Ka76}, we have
\begin{equation}\label{e:spectral_check}
\s(A)\cap N = \s\bigl( P_N(A)\circ A\circ P_N(A): \ran(P_N(A)) \to \ran(P_N(A))
\bigr).
\end{equation}

\begin{definition}\label{d:admissible} (Cf. Zhu \cite[Definition 1.3.6]{Zh00},
\cite[Definition 2.1]{Zh01}, and \cite[Definition 2.6]{ZhLo99}).
Let $\ell\<\C$ be a $C^1$ real 1-dimensional submanifold which has no boundary and
is co-oriented (i.e., with oriented normal bundle).
Let $A$ be a closed operator in a Banach space $X$.

\noi (a) We call $A$ {\em admissible} with respect to $\ell$, if
there exists a bounded open subset $N$ of $\C$ such that
\begin{equation}\label{e:a-admissible}
\s(A)\cap N=\s(A)\cap\ell,\  \s(A)\cap\partial N=\emptyset, \tand
\dim\ran P_N(A)<+\infty.
\end{equation}
Then $P_N(A)$ does not depend on the choice of such $N$. We set
\begin{equation}\label{e:nullity}
P_{\ell}(A):= P_N(A) \tand \nu_{\ell}(A):= \dim\ran P_N(A).
\end{equation}
For fixed $\ell$ and $X$ we shall denote the space of all $\ell$-admissible
closed operators in $X$ by $\Aa_{\ell}(X)$.

\smallskip

\noi (b) Let $A\in \Aa_{\ell}(X)$. Let $N\<\C$ be open and bounded with
$C^1$ boundary. We set ${N}^0={N}\cap \ell$ and assume
\begin{equation}\label{e:n-admissible}
\ol{{N}^0}=\ol {N}\cap \ell,\/
\s(A)\cap\ell\< N, \s(A)\cap\partial N=\emptyset, \tand \dim\ran P_N(A)<+\infty.
\end{equation}
Moreover, we require that each connected component of ${N}$ has connected intersection
with $\ell$ so that the disjoint positive (negative) part ${N}^\pm$ of ${N}$ with respect to
the co-orientation of $\ell$ is well-defined, and we have disjoint union $N=N^+ \cup N^0\cup   N^-$.
We shall call the resulting
triple $(N;N^+,N^-)$ {\em admissible} with respect to $\ell$ and $A$,
and write $(N;N^+,N^-)\in \Aa_{\ell,A}$\/.
\end{definition}

Now we are able to define spectral continuity and the spectral flow.
Our data are a co-oriented curve $\ell\<\C$, a family of Banach spaces
$\{X_s\}_{s\in [a,b]}$\, and a family $\{A_s\}_{s\in [a,b]}$ of
closed operators in $X_s$\,.

\begin{definition}\label{d:spectral_continuous}
(a) We shall call the family $\{A_s\}\in \Aa_{\ell}(X_s)$, $s\in [a,b]$
{\em spectral-continuous} near $\ell$ at $s_0\in [a,b]$, if
there is an $\e(s_0)>0$ such that for all $\e'\in(0, \e(s_0))$ there exists
a triple $(N;N^+,N^-)$ such that
\[
(N;N^+,N^-)\in\Aa_{\ell,A_s} \quad\text{ for all $\abs{s-s_0} <\e'$}\/;
\]
and for all triple
$(N';N'^+,N'^-)\in\Aa_{\ell,A_s}$ with $\ol{N'}\< N$, and $N'^\pm\< N^\pm$\/,
we have
\[
(N';N'^+,N'^-)\in\Aa_{\ell,A_s} \quad\text{ for all $\abs{s-s_0} \ll 1$}\/;
\]
and $\dim\range P_{N'}(A_s)$ and
$\dim\range P_{N^\pm\setminus N'^\pm}(A_s)$ do not depend on $s$.

\noi We shall call the family $\{A_s\}\in \Aa_{\ell}(X_s)$, $s\in [a,b]$
{\em spectral-continuous} near $\ell$, if it is spectral-continuous
near $\ell$ at $s_0$ for all $s_0\in [a,b]$.

\smallskip

\noi (b) Let $\{A_s\}$, $s\in [a,b]$ near $\ell$ be a spectral-continuous family.
Then there exist a partition
\begin{equation}\label{e:partition}
a=s_0\leq t_1\leq s_1\leq\ldots s_{n-1}\leq t_n\leq s_n=b
\end{equation}
of the interval $[a,b]$,
such that $s_{k-1},s_k\in (t_k-\e(t_k),t_k+\e(t_k))$, $k=1,\ldots,n$.
Let $(N_k;N_k^+,N_k^-)$
be like a $(N;N^+,N^-)$ in (a) for $t_k$ and some $\e'\in(0, \e(s_0))$ such that
$s_{k-1},s_k\in (t_k-\e',t_k+\e')$, $k=1,\ldots,n$.
Then we define the {\em spectral flow} of $\{A_s\}_{a\leq s\leq b}$ through $\ell$ by
\begin{multline}\label{e:spectral_flow}
\SF_{\ell}\bigl\{A_s; a\leq s \leq b\bigr\}
\\:= \sum_{k=1}^n\Bigl(\dim\ran\bigl(P_{N_k^-}(A_{s_{k-1}})\bigr)
- \dim\ran\bigl(P_{N_k^-}(A_{s_k})\bigr)\Bigr) .
\end{multline}

When $l=i\R$ with co-orientation from left to right, we set
$$\SF\bigl\{A_s; a\leq s \leq b\bigr\}
:=\SF_{\ell}\bigl\{A_s; a\leq s \leq b\bigr\}.$$

\end{definition}

Note that for a family of $\bigl\{A_s\in \Aa_{\ell}(X_s)\bigr\}$,
we always obtain a spectral-continuous family, when
we are given a suitable family of transformations $T_{s,s_0}:Y_s\to Y_{s_0}$
such that the family
\[
T_{s,s_0} A_s T_{s,s_0}\ii\in\Cc(Y_{s_0})
\]
is continuously varying.

From our assumptions
it follows that the spectral flow is independent of the choice of the
partition (\ref{e:partition}) and
admissible $(N_k;N_k^+,N_k^-)$, hence it is well defined.
From the definition it follows that the spectral flow through $\ell$
is path additive under catenation and homotopy invariant.
For details of the proof, see \cite{Ph96} and \cite{ZhLo99}.

\medskip

We close the appendix by discussing the invariance of the spectral flow
under embedding in a larger space, assuming a simple regularity condition.

\begin{lemma}\label{l:sf-embedding}
Let $\{Y_s; s\in [a,b]\}$ and $\{X_s; s\in [a,b]\}$ be two families of
(complex) Banach spaces with $X_s\< Y_s$ (no density or continuity of the
embeddings assumed). Let $\{A_s\in\Cc(Y_s); s\in [a,b]\}$ be a spectral-continuous
curve near a fixed co-oriented curve $\ell\<\C$. We assume that $A_s(X_s)\< X_s$
for all $s$ and that the curve $\{A_s|_{X_s}\in\Cc(Y_s); s\in [a,b]\}$ is also
spectral-continuous near $\ell$. Then we have
\[
\SF_{\ell}\{A_s;s\in [a,b]\}= \SF_{\ell}\{A_s|_{X_s};s\in [a,b]\}
\]
if the `regularity' $\nu_{\ell}(A_s)=\nu_{\ell}(A_s|_{X_s})$ holds for all $s\in
[a,b]$.
\end{lemma}

\smallskip

\begin{proof} We go back to the local definition of $\SF_{\ell}$ and reduce to the
finite-dimensional case. So, let $s_0\in [a,b]$. Choose a triple
\[
(N_1;N_1^+,N_1^-)\in \Aa_{\ell,A_{s_0}}
\]
such that $N_1$ satisfies
\eqref{e:a-admissible} for $A_{s_0}$\/. Then by spectral continuity,
there exists a triple $(N;N^+,N^-)$ with $\ol{N}\< N_1$ with
\[
(N;N^+,N^-)\in \Aa_{\ell,A_{s}} \quad\text{for $\abs{s-s_0}\ll 1$}.
\]
Then we have, again for $\abs{s-s_0}\ll 1$
\begin{equation}\label{e:nus}
\dim\range P_N(A_s)=\nu_{\ell}(A_{s_0})=\nu_{\ell}(A_s|_{X_{s_0}})
= \dim\range P_N(A_s|_{X_s})
\end{equation}
by spectral continuity and the regularity assumption.
Now we consider for each $\la\in\C\cap N$ the algebraic
multiplicities and find
\begin{equation}\label{e:kernels}
\dim \ker (A_s|_{X_s}-\la I|_{X_s})^k \leq \dim \ker (A_s-\la I)^k
\end{equation}
for each $k\in\N$. Comparing
\begin{align*}
\dim\range P_N(A_s)&=\sum_{\la\in\s(A_s)\cap N}\sum_{k\in\N} \dim\ker(A_s-\la I)^k\tand\\
\dim\range P_N(A_s|_{X_s})
&= \sum_{\la\in\s(A_s|_{X_s})\cap N}\sum_{k\in\N} \dim\ker(A_s|_{X_s}-\la I|_{X_s})^k
\end{align*}
we obtain from equation \eqref{e:nus} and the inequalities  \eqref{e:kernels} that
each term in the first and second preceding equation must coincide. So
\[
\s(A_s)\cap N=\s(A_s|_{X_s})\cap N;
\]
and the algebraic multiplicities with respect to $A_s$ and $A_s|_{X_s}$
coincide in each point.
By the definition of the spectral flow, the two spectral flows must coincide.
\end{proof}

\bigskip
\addtocontents{toc}{\medskip\noi}

\end{document}